\documentclass[11pt]{article}
\usepackage{preamble}

\title{Asymptotics for $t$-Core Partitions and Stanton's Conjecture}
\author{Matt Tyler}
\date{\vspace{-3.5ex}}

\begin{document}

\maketitle

\numberwithin{equation}{section}
\numberwithin{table}{section}
\numberwithin{figure}{section}
\numberwithin{thm}{section}

\begin{abstract}
A partition is a \textit{$t$-core partition} if $t$ is not one of its hook lengths. Let $c_t(N)$ be the number of $t$-core partitions of $N$. In 1999, Stanton conjectured $c_t(N) \le c_{t+1}(N)$ if $4 \le t \ne N-1$. This was proved for $t$ fixed and $N$ sufficiently large by Anderson, and for small values of $t$ by Kim and Rouse. In this paper, we prove Stanton's conjecture in general.

Our approach is to find a saddle point asymptotic formula for $c_t(N)$, valid in all ranges of $t$ and $N$. This includes the known asymptotic formulas for $c_t(N)$ as special cases, and shows that the behavior of $c_t(N)$ depends on how $t^2$ compares in size to $N$. For example, our formula implies that if $t^2 = \kappa N + o(t)$, then $c_t(N) = \frac{\exp\prn{2\pi\sqrt{A N}}}{B N} \prn{1 + o(1)}$ for suitable constants $A$ and $B$ defined in terms of $\kappa$.
\end{abstract}

\section{Introduction and statement of results}
Let $\l = (\l_1, \dots, \l_k)$ with $\l_1 \ge \dots \ge \l_k > 0$ be a partition of $|\l| = \l_1 + \dots + \l_k$. The \textit{Ferrers--Young diagram} of $\l$ consists of $k$ rows of left-justified boxes, with $\l_i$ boxes in row $i$. Let $\l_j'$ be the number of boxes in column $j$ (so $\l' = (\l_1', \dots, \l_\ell')$ is the conjugate partition to $\l$). The \textit{hook length} for the box at row $i$ and column $j$ is $(\l_i - j) + (\l_j' - i) + 1$. In other words, the hook length is the number of boxes directly to the right, plus the number of boxes directly below, plus $1$ for the box itself. We say $\l$ is a \textit{$t$-core partition} if $t$ is not one of its hook lengths. See Figure \ref{fig:young-diagram} for an example.

\begin{figure}[h]
\captionsetup{width=0.8\textwidth}
\[
\begin{ytableau}
8 & 7 & 5 & 4 & 2 & 1 \\
5 & 4 & 2 & 1 \\
2 & 1
\end{ytableau}
\]
\caption{The Ferrers--Young diagram for the partition $(6, 4, 2)$ of $12$, labelled with hook lengths. This is a $t$-core partition for $t \not\in \{1, 2, 4, 5, 7, 8\}$.}
\label{fig:young-diagram}
\end{figure}

We write $c_t(N)$ for the number of $t$-core partitions of $N$, and $p(N)$ for the number of partitions of $N$.

Much of the motivation for studying $t$-core partitions comes from the representation theory of the symmetric group $S_N$. The irreducible complex representations of $S_N$ are parameterized by the partitions $\l$ of $N$. The representation associated to $\l$ has dimension $\frac{N!}{\prod_{h\in\mathcal{H}(\l)} h}$ where $\mathcal{H}(\l)$ is the multiset of hook lengths of $\l$.

When $p$ is prime, a $p$-core partition of $N$ is therefore a partition for which the power of $p$ dividing $|S_N|$ is the same as the power of $p$ dividing the dimension of the associated representation. These are the projective irreducible representations in characteristic $p$, i.e. the $p$-blocks of defect $0$. More generally, by Nakayama's conjecture (proved by Brauer and Robinson), $p$-core partitions parameterize the $p$-blocks of representations of $S_n$. See \cite[sections 2.7 and 6.1]{james-kerber} for details.

Another application of $t$-core partitions is due to Garvan, Kim, and Stanton \cite{garvan90}, who used $t$-core partitions for $t = 5, 7, 11$ to find combinatorial proofs of Ramanujan's congruences
\aleq
p(5N+4) &\equiv 0 \Mod{5}, \\
p(7N+5) &\equiv 0 \Mod{7}, \quad \text{and} \\
p(11N+6) &\equiv 0 \Mod{11}.
\alqe

The \textit{$t$-core partition conjecture} states that $c_t(N) > 0$ for all $t \ge 4$ and $N \ge 1$. This was proved for $t = 5, 7$ by Erdmann and Michler \cite{erdmann}, when $t$ has a prime factor $\le 11$ by Ono \cite{ono94, ono95}, and in general by Granville and Ono \cite{granville00}. In 1999, Stanton \cite{stanton} made the ``possibly rash" conjecture that in fact much more is true.
\begin{con}[Stanton] \label{conj:stanton}
For all $4 \le t \ne N - 1$, $c_t(N) \le c_{t+1}(N)$.
\end{con}

We always have $c_t(N) \le p(N)$, with equality if $t > N$, so the conjecture is obvious if $t \ge N$. The restriction $t \ne N-1$ is necessary because $c_{N-1}(N) = p(N) - N + 1 = c_N(N) + 1$.

In the case $t > \frac{N}{2}$, Stanton's conjecture was quickly proved by Craven \cite{craven} using combinatorial arguments. In 1999, Lulov and Pittel \cite{lulov} used the circle method to show that if $t \ge 6$, then
\eq \label{eqn:anderson}
c_t(N)
= A_t(N) \frac{(2\pi)^{\frac{t-1}{2}}}{t^{\frac{t}{2}} \G\prn{\frac{t-1}{2}}} \prn{N + \frac{t^2-1}{24}}^{\frac{t-3}{2}} + O_t\prn{N^{\frac{t-1}{4}}}
\qe
where $A_t(N)$ is an arithmetic quantity satisfying $\abs{A_t(N) - 1} \le \abs{\z(\frac{t-3}{2})-1}$. This was rediscovered in 2008 by Anderson \cite{anderson}, who used it to prove the conjecture when $t$ is fixed and $N$ is sufficiently large. This work was made explicit by Kim and Rouse \cite{kim-rouse} in 2014, who proved the conjecture for certain values of $t$.\footnote[2]{There may be a technical error in their proof, but it works without modification for $t \ll \log N$, and may be extendable up to $t \ll N^{\frac{1}{4}}$. For completeness, we will independently treat the case $t \ge 8$ here, and rely on their work for $4 \le t < 8$, which is unaffected by the error.} 

In this paper, we will prove Stanton's conjecture.
\begin{thm}
Stanton's conjecture \ref{conj:stanton} is true.
\end{thm}
\begin{RMK}
Our proof shows that the inequality is strict for $4 \le t < N - 1$, except for $c_5(10) = c_6(10) = 12$.
\end{RMK}

Previously, $c_t(N)$ has only been understood when $t$ is of constant size, as in (\ref{eqn:anderson}), or when $t > (1+\e) \frac{\sqrt{6}}{2\pi} \sqrt{N} \log N$, in which case $c_t(N) \sim p(N)$ \cite{morotti}. Our contribution is a saddle point asymptotic formula for $c_t(N)$ that is valid in all ranges of $t$ and $N$ and includes these formulas as special cases. Before stating the theorem, we will need a few definitions.

Let $e(z) = e^{2 \pi i z}$, and for $\Im z > 0$, let $\eta(z)$ denote the Dedekind eta function
\eq \label{eqn:dedekind-eta}
\eta(z) = e\prn{\frac{z}{24}} \prod_{n=1}^\infty (1 - e(n z)).
\qe
Let
\eq
f_t(z) = \frac{\eta(tz)^t}{\eta(z)},
\qe
which is a modular form of weight $\frac{t-1}{2}$ and level $t$, and which will be of central importance in this paper. The relationship between $f_t(z)$ and $t$-core partitions comes from the remarkable classical identity for the generating function of $c_t(N)$ \cite[(2.7.17)]{james-kerber}
\eq
\sum_{N=0}^\infty c_t(N) e(Nz)
= \prod_{n=1}^\infty \frac{(1 - e(n t z))^t}{1 - e(nz)}
= e\prn{\tfrac{1-t^2}{24} z} f_t(z).
\qe

We also define
\eq \label{eqn:func-def}
\func_k(z) = - \frac{z^{k+1}}{2 \pi i} \prn{\dv{}{z}}^k \log \eta(z).
\qe
When $k \ge 1$, $\func_k(z) = -\frac{z^{k+1}}{24} E_2^{(k-1)}(z)$ where $E_2(z) = 1 - 24 \sum_{n=1}^\infty \s(n) e(nz)$ is the Eisenstein series of weight $2$ and $\s(n) = \sum_{d\divides n} d$. In particular, note that $\func_1(iy) = \frac{y^2}{24} \prn{1 - 24 \sum_{n=1}^\infty \s(n) e^{-2\pi n y}}$ is real-valued.

Our main theorem is the following asymptotic formula for $c_t(N)$. Unless otherwise indicated with a subscript, all implied constants are universal.

\begin{thm} \label{thm:t-core-asymptotic}
Suppose $N \ge 0$ and $t \ge 6$.
\begin{enumerate}[(i)]
    \item \label{item:y} There is a unique value of $y > 0$ satisfying
    \[\frac{\func_1(ity) - \func_1(iy)}{y^2} = N + \frac{t^2-1}{24}.\]
    Moreover, $y$ lies in the range
    \[
    \frac{t-1}{4 \pi \prn{N + \frac{t^2-1}{24}}}
    < y <
    \frac{1}{\frac{3}{\pi} + \sqrt{24 N - 1 + \frac{9}{\pi^2}}}.
    \]
    \item \label{item:asymp} For this value of $y$,
    \[
    c_t(N) = \frac{y^{\frac{3}{2}} e^{2 \pi y \prn{N + \frac{t^2-1}{24}}} f_t(iy)}{\sqrt{\func_2(iy) - \func_2(ity)}}
    \prn{
    1 + O\prn{\frac{1}{\min(t, \sqrt{N})}}
    }.
    \]
    \item \label{item:diff-asymp} For this value of $y$, if $t y \ge \frac{1}{2}$, then
    \[
    c_{t+1}(N+t) - c_t(N+t)
    = (2\pi t y - 1) c_t(N) \prn{1 + O\prn{
    y + t y e^{-2\pi t y}
    }}.
    \]
\end{enumerate}
\end{thm}

Theorem~\ref{thm:t-core-asymptotic} shows that $c_t(N)$ exhibits a change in behavior depending on how $t^2$ compares in size to $N$. In order to interpret this result, let us examine in turn the ranges where $t^2$ is much smaller than $N$, similar in size to $N$, and much bigger than $N$.

In the range $t^2 \log t < \frac{8 \pi^2}{1+\e} N$ where $\e > 0$, we have $y \approx \frac{t-1}{4\pi \prn{N + \frac{t^2-1}{24}}}$, and part \ref{item:asymp} of Theorem~\ref{thm:t-core-asymptotic} implies
\eq \label{eqn:small-t-corollary}
c_t(N)
= \frac{(2\pi)^{\frac{t-1}{2}}}{t^{\frac{t}{2}} \G\prn{\frac{t-1}{2}}} \prn{N + \frac{t^2-1}{24}}^{\frac{t-3}{2}} \prn{1 + O\prn{t^{-\min\prn{\e, 1}}}}.
\qe
This expands the range of validity of Lulov and Pittel's formula (\ref{eqn:anderson}) from fixed $t$ to $t^2 \log t < \frac{8 \pi^2}{1+\e} N$. When $t^2 \log t < \frac{8 \pi^2}{2+\e} N$, a slightly different approach allows us to refine the error term in (\ref{eqn:small-t-corollary}) to be exponentially small.

\begin{thm} \label{thm:t-core-asymptotic-small-t}
For all $\e > 0$, if $t \ge 6$ and $t^2 \log t < \frac{8\pi^2}{2+\e} N$, then
\[
c_t(N)
= \frac{(2\pi)^{\frac{t-1}{2}}}{t^{\frac{t}{2}} \G\prn{\frac{t-1}{2}}} \prn{N + \frac{t^2-1}{24}}^{\frac{t-3}{2}}
\prn{1 + O\prn{
    t^{-\frac{1}{2}} (1+\e)^{-\frac{1}{4} t}
    + e^{-\frac{1}{6} t}
    + t e^{-\frac{8\pi^2}{t(t-1)} N}
    }
}.
\]
\end{thm}

In the range $t > (1+\e) \frac{\sqrt{6}}{2\pi} \sqrt{N} \log N$, we have $y \approx \frac{1}{\frac{3}{\pi} + \sqrt{24 N - 1 + \frac{9}{\pi^2}}}$, and part \ref{item:asymp} of Theorem~\ref{thm:t-core-asymptotic} tells us only that
\eq \label{eqn:t-core-asymptotic-big-t}
c_t(N) = p(N) \prn{1 + O\prn{N^{-\frac{1}{2}} + N^{-\frac{\e}{2}} \log N}}.
\qe
In this range, almost all partitions are $t$-cores, and we may once again refine (\ref{eqn:t-core-asymptotic-big-t}), this time to evaluate the remainder.

\begin{thm} \label{thm:t-core-asymptotic-big-t}
For all $\e > 0$, if $t > (1+\e) \frac{\sqrt{6}}{2\pi} \sqrt{N} \log N$, then
\eq
p(N) - c_t(N) = t p(N-t) + O_\e\prn{t^2 p(N-2t)}.
\qe
\end{thm}

In the intermediary range $t \asymp \sqrt{N}$, $t$-core partitions have not previously been understood. It turns out that $c_t(N)$ has an asymptotic formula of the same shape as the Hardy--Ramanujan asymptotic formula for partitions.

\begin{cly} \label{cly:t-core-asymptotic-middle-t}
Suppose $t^2 = \kappa N + o(t)$. Then,
\[
c_t(N) = \frac{e^{2 \pi \sqrt{A(\kappa) N}}}{B(\kappa) N} \prn{1 + o(1)}
\]
where $A(\kappa)$ and $B(\kappa)$ are constants defined by choosing $v$ so that
\[
\frac{1}{24 v^2} - \frac{1}{24} + \frac{1}{v^2} \func_1(iv) = \frac{1}{\kappa},
\]
and letting
\[
A(\kappa) = \frac{\kappa}{v^2} \prn{\frac{1}{12} - \func_0(iv) + \func_1(iv)}^2 \quad \text{and} \quad
B(\kappa) = \frac{\kappa}{v^2} \sqrt{\frac{1}{12} - \func_2(iv)}.
\]
\end{cly}
We will show in \S \ref{ssec:middle-t-proof} that the function $\frac{1}{24 v^2} - \frac{1}{24} + \frac{1}{v^2} \func_1(iv)$ is a monotonically decreasing bijection $(0, \infty) \to (0, \infty)$, so $v$ always exists and is unique.

\begin{figure}[h]
  \centering
  \footnotesize
  \begin{minipage}[b]{0.45\textwidth}
    \captionsetup{width=0.8\textwidth, skip=0pt}
    \begin{tikzpicture}
        \node at (0,0)
            {\includegraphics[width=6cm]{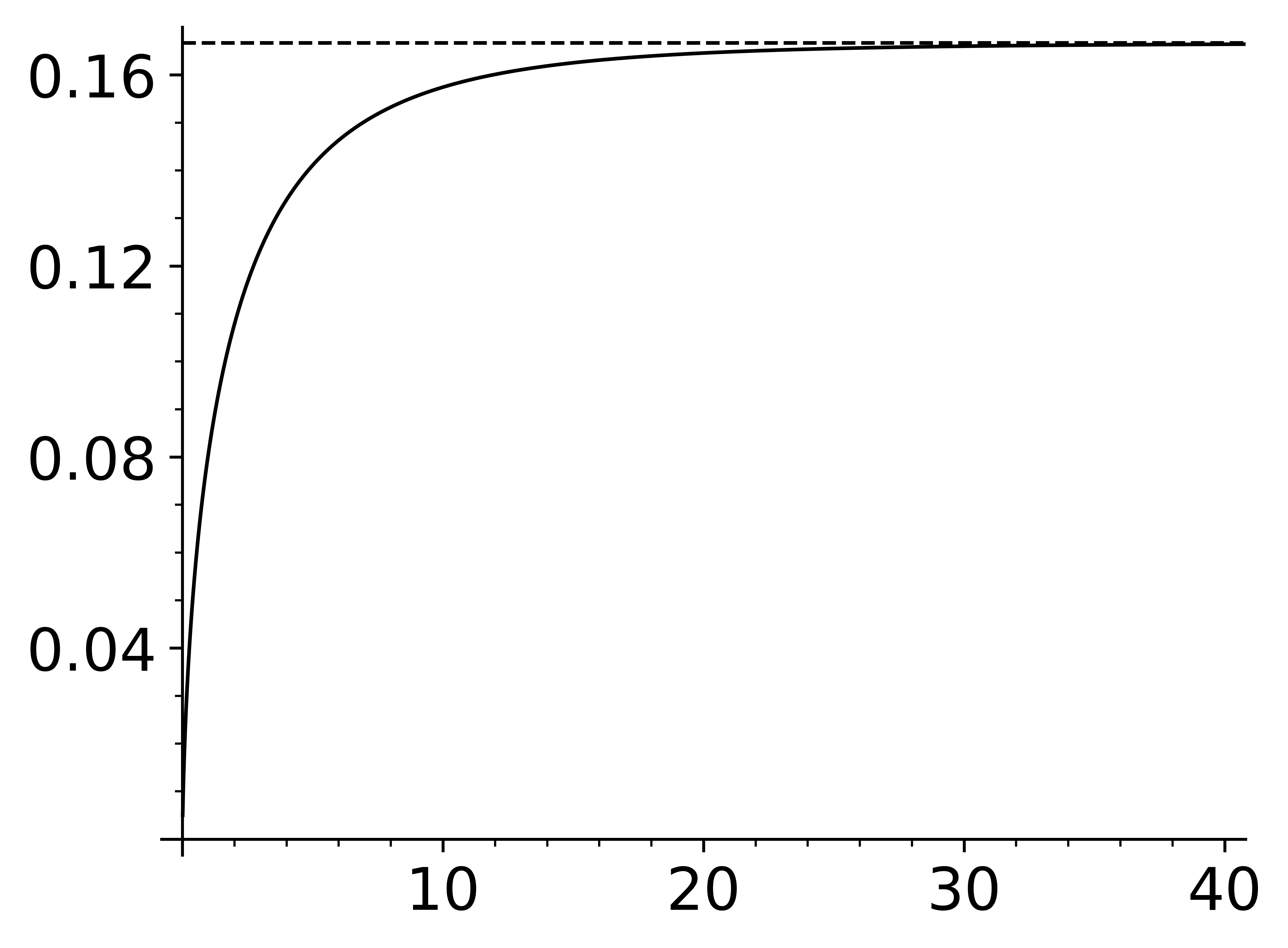}};
        \fill [white] (-2.24,-1.2) rectangle (-2.95, 2.05);
        \node at (-2.57,-.81) {0.04};
        \node at (-2.57,0.08) {0.08};
        \node at (-2.57,.97) {0.12};
        \node at (-2.57,1.86) {0.16};
        \fill [white] (-1.2,-2.08) rectangle (2.9,-1.78);
        \node at (-0.93,-1.93) {10};
        \node at (0.28,-1.93) {20};
        \node at (1.49,-1.93) {30};
        \node at (2.7,-1.93) {40};
    \end{tikzpicture}
    \caption{A graph of $A(\kappa)$ with an asymptote at $\frac{1}{6}$.}
  \end{minipage}
  \hfill
  \begin{minipage}[b]{0.45\textwidth}
    \captionsetup{width=0.8\textwidth, skip=0pt}
    \begin{tikzpicture}
        \node at (0,0)
            {\includegraphics[width=6cm]{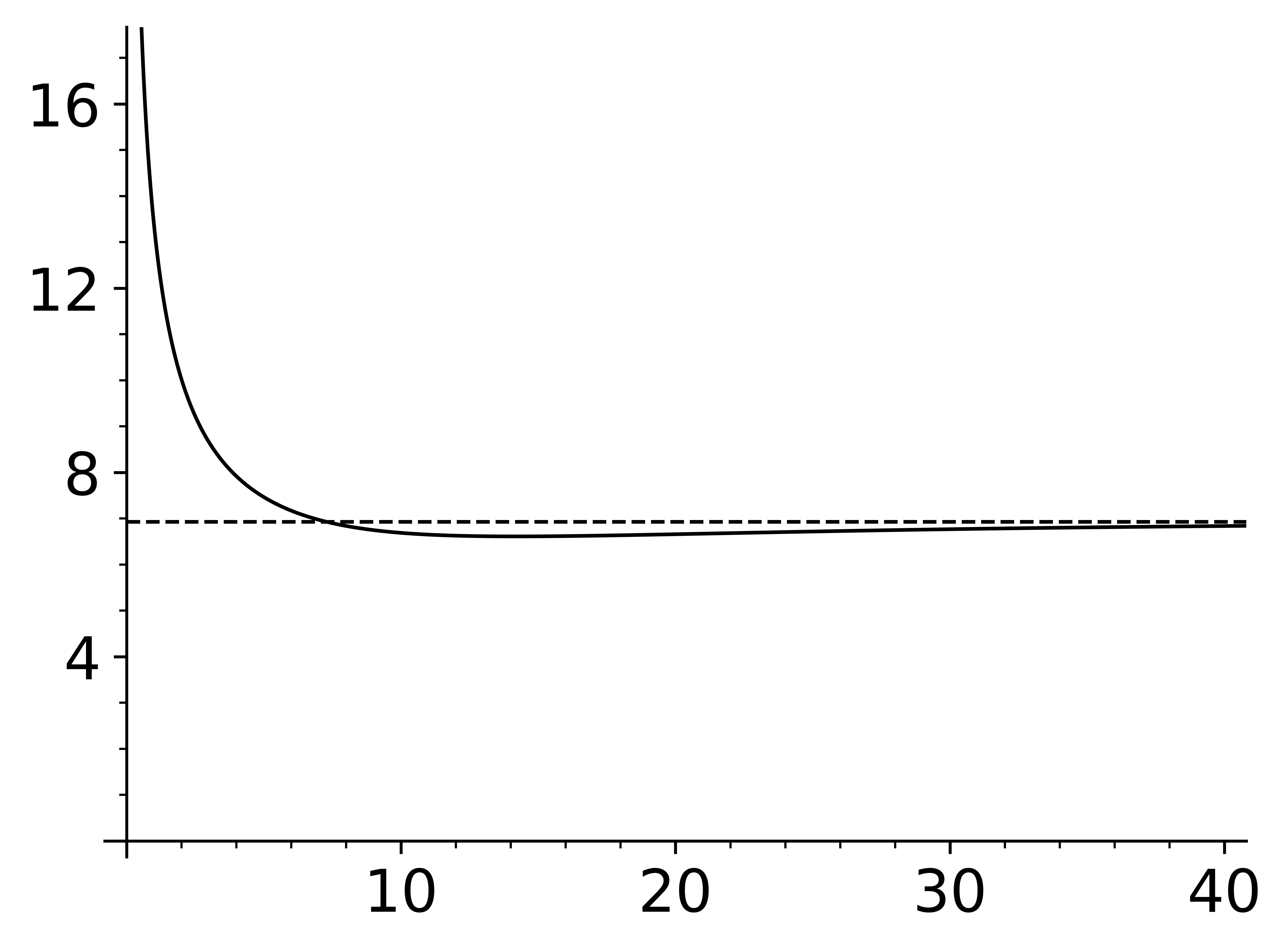}};
        \fill [white] (-2.49,-1.1) rectangle (-2.9, 2);
        \node at (-2.61,-.875) {4};
        \node at (-2.61,0) {8};
        \node at (-2.7,.875) {12};
        \node at (-2.7,1.75) {16};
        \fill [white] (-1.35,-2.08) rectangle (2.9,-1.79);
        \node at (-1.14,-1.93) {10};
        \node at (0.14,-1.93) {20};
        \node at (1.42,-1.93) {30};
        \node at (2.7,-1.93) {40};
    \end{tikzpicture}
    \caption{A graph of $B(\kappa)$ with an asymptote at $4\sqrt{3}$.}
  \end{minipage}
\end{figure}

The constants $A(\kappa)$ and $B(\kappa)$ of Corollary \ref{cly:t-core-asymptotic-middle-t} satisfy
\eq
\lim_{\kappa\to\infty} A(\kappa) = \frac{1}{6}
\quad \text{and} \quad
\lim_{\kappa\to\infty} B(\kappa) = 4 \sqrt{3},
\qe
which makes sense in light of (\ref{eqn:t-core-asymptotic-big-t}) because $p(N) = \frac{e^{2\pi \sqrt{N/6}}}{4\sqrt{3} N} \prn{1 + o(1)}$.

When $\frac{t^2}{N}$ is bounded from above, the ratio of main terms for $c_{t+1}(N)$ and $c_t(N)$ in part \ref{item:asymp} of Theorem~\ref{thm:t-core-asymptotic} is uniformly bounded from below by a constant larger than $1$. This will prove Stanton's conjecture in this range for sufficiently large $t$. When $\frac{t^2}{N}$ grows larger, (\ref{eqn:t-core-asymptotic-big-t}) presents an obstruction --- both $c_{t+1}(N)$ and $c_t(N)$ are approximately the same size as $p(N)$. In this case, it is better to consider their difference $c_{t+1}(N) - c_t(N)$ rather than their ratio, and we will use part \ref{item:diff-asymp} of Theorem~\ref{thm:t-core-asymptotic} to complete the proof of Stanton's conjecture for $t$ sufficiently large.

Since our goal is to prove Stanton's conjecture for all $t$, not just $t$ sufficiently large, we will actually keep track of the constants involved and prove an explicit version of Theorem~\ref{thm:t-core-asymptotic} in Theorems~\ref{thm:t-core-asymptotic-explicit} and \ref{thm:t-core-diff-asymptotic-explicit} and an explicit version of Theorem~\ref{thm:t-core-asymptotic-small-t} in Remark~\ref{rmk:t-core-asymptotic-small-t-explicit}. This will allow us to precisely determine which values of $t$ and $N$ are large enough for the asymptotic formulas to settle the conjecture. The remaining cases will be few enough that we can check them manually with SageMath \cite{sagemath}.

The organization of this paper is as follows. In \S \ref{sec:sketch}, we will sketch the argument used to prove Theorem \ref{thm:t-core-asymptotic} and motivate why the main terms take the form they do. In \S \ref{sec:away-from-0}, we will bound the contribution away from the saddle point. This is the heart of the paper, and is the main input that allows us to evaluate $c_t(N)$ in new ranges. In \S \ref{sec:technical-lemmas}, we will prove some technical lemmas that will be useful in \S \ref{sec:asymptotic-formulas}, mostly involving bounding the functions $\func_k$ for $k \le 4$. In \S \ref{sec:asymptotic-formulas}, we will use the saddle point method to prove explicit versions of Theorems~\ref{thm:t-core-asymptotic} and \ref{thm:t-core-asymptotic-small-t} and also prove Theorem~\ref{thm:t-core-asymptotic-big-t} and Corollary \ref{cly:t-core-asymptotic-middle-t}. In \S \ref{sec:stanton}, we will prove Stanton's conjecture.

\section*{Notation}

Because our goal is to prove an exact inequality, we will require the use of an explicit version of big $O$ notation. To that end, we use the symbol $\err$ to denote a complex number of magnitude $\le 1$ which is allowed to depend on all parameters, and which may be different at each occurrence. For example, $a = b + \err c$ means $|a - b| \le |c|$. As usual with big $O$ notation, equations involving $\err$ should be read from left to right, so for example $1 = 2 \err$, but it would not be correct to write $2 \err = 1$.

We will always let $z = x + i y$ where $x$ and $y$ are real numbers and $y > 0$.

We also make use of the shorthand
\eq
M = N + \frac{t^2-1}{24}.
\qe

\section*{Acknowledgments}

I am grateful to my advisor Kannan Soundararajan for suggesting this problem and for offering invaluable advice. I would also like to thank Shintaro Fushida-Hardy, Byungchan Kim, Milo Marsden, Jeremy Rouse, Dennis Stanton, Zachary Stier, and the two anonymous referees for many helpful comments. This paper is based upon work supported by the National Science Foundation Graduate Research Fellowship under Grant No. DGE-1656518.

\section{Sketch of the proof of Theorem \ref{thm:t-core-asymptotic}}
\label{sec:sketch}
Let us begin by sketching the method used to prove Theorem \ref{thm:t-core-asymptotic}. Over the course of \S \ref{sec:away-from-0}, \S \ref{sec:technical-lemmas}, and \S \ref{sec:asymptotic-formulas}, we will make this argument rigorous and explicit.

In contrast to previous papers on this subject \cite{granville00, anderson, kim-rouse}, which decompose $f_t(z)$ into a sum of Eisenstein series and cusp forms, our approach is to use the saddle point method to evaluate the integral
\eq \label{eqn:t-core-integral-intro}
c_t(N) = \int_{-1/2}^{1/2} e\prn{-M z} f_t(z) dx = \int_{-1/2}^{1/2} e\prn{-M z + \frac{1}{2\pi i} \log f_t(z)} dx.
\qe

The functions $\func_k(z)$ defined in (\ref{eqn:func-def}) arise from the term $\frac{1}{2\pi i} \log f_t(z)$ in the exponential above, whose Taylor expansion is
\eq \label{eqn:func-taylor-expansion}
\frac{1}{2\pi i} \log f_t(z)
= \sum_{k=0}^\infty \frac{x^k}{k!} \frac{\func_k(iy) - \func_k(ity)}{(i y)^{k+1}}.
\qe

In \S \ref{sec:away-from-0}, we will use the functional equation for $\eta(z)$ to show that this integral is dominated by the contribution when $x$ is close to $0$, and in \S \ref{sec:technical-lemmas}, we will show that the higher order terms in the expansion (\ref{eqn:func-taylor-expansion}) are small compared to the second order term.

Together, these estimates imply that we may approximate $c_t(N)$ as
\eq \label{eqn:simplified-t-core-integral}
c_t(N) \approx \int_{-y/3}^{y/3} e\prn{
    -M (x+iy)
    + \sum_{k=0}^2 \frac{x^k}{k!} \frac{\func_k(iy) - \func_k(ity)}{(i y)^{k+1}}
} dx.
\qe
So far, this approximation holds for any $y > 0$, and we are free to choose whichever value of $y$ makes (\ref{eqn:simplified-t-core-integral}) easiest to evaluate. To that end, note that the integrand in (\ref{eqn:simplified-t-core-integral}) is an exponential of a quadratic polynomial in $x$, and if we choose $y$ satisfying $\frac{\func_t(ity) - \func_1(iy)}{y^2} = M$, as in part \ref{item:y} of Theorem~\ref{thm:t-core-asymptotic}, then the linear term vanishes. With this value of $y$, our integral reduces to
\aleq
c_t(N)
\approx e^{2\pi M y} f_t(iy) \int_{-y/3}^{y/3} e\prn{
\frac{x^2}{2} \frac{\func_2(iy) - \func_2(ity)}{(iy)^3}
} dx.
\alqe

If we expand the range of integration now to $(-\infty, \infty)$, then what remains is simply a Gaussian integral, which we can evaluate to find
\eq
c_t(N)
\approx e^{2\pi M y} f_t(iy) \sqrt{\frac{y^3}{\func_2(iy) - \func_2(ity)}}.
\qe
This is the main term in part \ref{item:asymp} of Theorem~\ref{thm:t-core-asymptotic}.

For part \ref{item:diff-asymp} of Theorem~\ref{thm:t-core-asymptotic} we must analyze the difference
\eq
c_{t+1}(N+t) - c_t(N+t)
= \int_{-1/2}^{1/2} e\prn{-(M + t) z} f_t(z) \prn{e\prn{-\frac{2t+1}{24} z} \frac{f_{t+1}(z)}{f_t(z)} - 1} dx.
\qe

Let us consider the last term $e\prn{-\frac{2t+1}{24} z} \frac{f_{t+1}(z)}{f_t(z)} - 1$ in the integrand. So long as $z$ is small enough (once again afforded by \S \ref{sec:away-from-0}) and $t y$ is large enough, we can use the approximation
\aleq
e\prn{-\frac{2t+1}{24} z} \frac{f_{t+1}(z)}{f_t(z)}
&= \prod_{n=1}^\infty \frac{\prn{1 - e(n(t+1) z)}^{t+1}}{\prn{1 - e(n t z)}^t} \\
&\approx 1 + t e(t z) - (t+1) e((t+1) z) \\
&\approx 1 - e(t z) \prn{2 \pi i t z + 1}.
\alqe
Therefore,
\aleq
c_{t+1}(N+t) - c_t(N+t)
&\approx \int_{-1/2}^{1/2} e\prn{-(M + t) z} f_t(z) (-2\pi i t z - 1) e(tz) dx \\
&= \int_{-1/2}^{1/2} e\prn{-M z} f_t(z) (2\pi t y - 1 - 2\pi i t x) dx \\
&= (2 \pi t y - 1) c_t(N) - 2 \pi i t \int_{-1/2}^{1/2} x e\prn{-M z} f_t(z) dx.
\alqe
Since our chosen value of $y$ makes $e(-Mz) f_t(z)$ approximately an even function of $x$, this last integral is small, and we obtain
\eq
c_{t+1}(N+t) - c_t(N+t)
\approx (2 \pi t y - 1) c_t(N).
\qe
This is the main term in part \ref{item:diff-asymp} of Theorem~\ref{thm:t-core-asymptotic}.

\section{Bounding the contribution away from the saddle point}
\label{sec:away-from-0}
Our starting point is the Fourier inversion integral
\eq \label{eqn:t-core-integral-intro-2}
c_t(N) = \int_{-1/2}^{1/2} e\prn{-M z} f_t(z) dx.
\qe
The Fourier coefficients $c_t(N)$ of $f_t(z)$ are all non-negative, so $\abs{f_t(z)} \le \abs{f_t(iy)} = f_t(iy)$, and hence the magnitude of the integrand in (\ref{eqn:t-core-integral-intro-2}) is maximized when $x = 0$. Our goal in this section is to prove the following proposition, which shows that this estimate can be substantially improved so long as $x$ is not too close to $0$.


\begin{prop} \label{prop:minor-arc-bound}
If $\min(t, \frac{1}{y}) \ge 1000$, then
\[
\int_{\frac{y}{3} \le \abs{x} \le \frac{1}{2}} \abs{\frac{f_t(z)}{f_t(iy)}} dx
< y e^{-\frac{1}{70} \min(t, \frac{1}{y})}.
\]
\end{prop}

Recall that the group $\SL_2(\Z)$ of $2\times2$ integer matrices with determinant $1$ acts on the upper half plane $\H = \{z \mid y > 0\}$ by linear fractional transformations
\eq
\g z = \frac{a z + b}{c z + d} \text{ for } \g = \Matrix{a & b \\ c & d} \in \SL_2(\Z).
\qe
The Dedekind eta function (\ref{eqn:dedekind-eta}) satisfies the transformation formula
\eq
\eta\prn{\frac{az+b}{cz+d}} = e\prn{\frac{a+d}{24c} + \half s(-d, c)} \sqrt{-i (cz+d)} \eta(z)
\text{ for all } \Matrix{a & b \\ c & d} \in \SL_2(\Z)
\qe
where $s(h, k)$ is the Dedekind sum $\sum_{r=1}^{k-1} \frac{r}{k} \prn{\set{\frac{hr}{k}} - \half}$. This enables us to approximate $\eta(z)$ anywhere in the upper half plane, and in particular implies
\aleq \label{eqn:eta-feq}
\eta(-1/z) &= \sqrt{-iz} \eta(z)
\quad \text{and} \\
\abs{\eta(\g z)} &= \prn{\frac{y}{\Im \g z}}^{\frac{1}{4}} \abs{\eta(z)} \text{ for all } \g \in \SL_2(\Z)
\alqe
by applying the matrix $\sMatrix{0 & -1 \\ 1 & 0}$ and taking magnitudes, respectively.

We will need the following two lemmas for the proof of Proposition \ref{prop:minor-arc-bound}.

\begin{lma} \label{lma:eta-bound-1}
For all $x \in \R$ and $y \ge \frac{\sqrt{3}}{2}$,
\[
\abs{\eta(z)} = e^{-\frac{\pi}{12} y + \err 1.01 e^{-2\pi y}}.
\]
\end{lma}

\begin{proof}
Suppose $y \ge \frac{\sqrt{3}}{2}$. Since $\log \eta(z) = \frac{\pi i}{12} z - \sum_{n=1}^\infty \frac{\s(n)}{n} e(nz)$, we may use the elementary bound $\s(n) \le n^2$ to find
\eq
\abs{\log \eta(z) - \frac{\pi i}{12} z}
\le \sum_{n=1}^\infty \frac{\s(n)}{n} e^{-2 \pi y n}
< \sum_{n=1}^\infty n e^{-2 \pi y n}
\le \frac{e^{-2\pi y}}{(1 - e^{-\sqrt{3} \pi})^2}
< 1.01 e^{-2\pi y}.
\qe
\end{proof}

\begin{lma} \label{lma:eta-bound-2}
For all $x \in \R$ and $y > 0$,
\[
\abs{\eta(z)} < \frac{7}{9} y^{-\frac{1}{4}}.
\]
\end{lma}

In fact, $y^{\frac{1}{4}} \abs{\eta(z)}$ is maximized when $z = \frac{1+i\sqrt{3}}{2}$, where it equals $\prn{\frac{\sqrt{2} \pi}{3 \Gamma(2/3)^3}}^{\frac{1}{2}} = 0.7723\dots$ \cite{rouse08}.

\begin{proof}
By the second line in (\ref{eqn:eta-feq}), the function $(\Im z)^{\frac{1}{4}} \abs{\eta(z)}$ is $\SL_2(\Z)$-invariant, so we may assume $z = x + iy$ is in the standard fundamental domain $\set{z \in \H \mid \abs{z} \ge 1, \abs{x} \le \half}$ for $\SL_2(\Z)$, and in particular $y \ge \frac{\sqrt{3}}{2}$. By Lemma \ref{lma:eta-bound-1},
\eq
y^{\frac{1}{4}} \abs{\eta(x+iy)} \le y^\frac{1}{4} e^{-\frac{\pi}{12} y + 1.01 e^{-\sqrt{3} \pi}} \le \prn{\frac{3}{\pi}}^{\frac{1}{4}} e^{-\frac{\pi}{12} \frac{3}{\pi} + 1.01 e^{-\sqrt{3} \pi}} < \frac{7}{9}.
\qe
\end{proof}

The idea behind the proof of Proposition~\ref{prop:minor-arc-bound} is that Lemma~\ref{lma:eta-bound-1} provides us with a good bound for $\eta(z)$ so long as $y$ is large, so we choose $\g \in \SL_2(\Z)$ so that $\Im \g z$ is large, and use (\ref{eqn:eta-feq}) to relate $\abs{\eta(\g z)}$ with $\abs{\eta(z)}$. Sometimes, we are lucky and $\Im \g t z$ is also large. If not, we will need to find another element of $\SL_2(\Z)$ whose action on $t z$ sends it to a value with large imaginary part.

\begin{proof}[Proof of Proposition \ref{prop:minor-arc-bound}]
Choose $\g = \sMatrix{a & b \\ c & d} \in \SL_2(\Z)$ such that $\Im \g z \ge \frac{\sqrt{3}}{2}$, and assume without loss of generality that $c > 0$ and $0 \le d < c$. By (\ref{eqn:eta-feq}) and Lemma \ref{lma:eta-bound-1},
\aleq
\eta(iy) &= y^{-\frac{1}{2}} \eta(i/y) \le y^{-\frac{1}{2}} e^{-\frac{\pi}{12} \frac{1}{y} + 1.01 e^{-2000\pi}}
\quad \text{and} \\
\abs{\eta(z)}
&= \prn{\frac{y}{\Im \g z}}^{-\frac{1}{4}} \abs{\eta(\g z)}
\ge \prn{\frac{y}{\Im \g z}}^{-\frac{1}{4}} e^{-\frac{\pi}{12} \Im \g z - 1.01 e^{-\sqrt{3} \pi}},
\alqe
so we obtain
\eq
\abs{\frac{\eta(iy)}{\eta(z)}} < \prn{y \Im \g z}^{-\frac{1}{4}} e^{-\frac{\pi}{12} \prn{\frac{1}{y} - \Im \g z} + \frac{1}{228}}.
\qe

It remains to bound $\abs{\frac{\eta(tz)^t}{\eta(ity)^t}}$, for which we split up the proof into two cases depending on whether $ty \ge 1$ or $ty \le 1$.

First, suppose $ty \ge 1$, so that by Lemma \ref{lma:eta-bound-1},
\eq
\abs{\frac{\eta(tz)^t}{\eta(ity)^t}}
= e^{\err 2.02 t e^{-2\pi t y}}
\le e^{2.02 e^{-2\pi} \frac{1}{y}}
< e^{\frac{1}{265} \frac{1}{y}},
\qe
and hence
\eq
\abs{\frac{f_t(z)}{f_t(iy)}}
= \abs{\frac{\eta(tz)^t}{\eta(ity)^t}} \abs{\frac{\eta(iy)}{\eta(z)}}
< \prn{y \Im \g z}^{-\frac{1}{4}} e^{-\frac{\pi}{12} \prn{\frac{1}{y} - \Im \g z - \frac{12}{\pi} \frac{1}{265} \frac{1}{y}} + \frac{1}{228}}
\qe

Note that $y \Im \g z = \frac{y^2}{(c x + d)^2 + c^2 y^2}$. If $c > 1$, then $y \Im \g z \le \frac{1}{c^2} \le \frac{1}{4}$. Otherwise, $c = 1$ and $d = 0$, and $y \Im \g z = \frac{y^2}{x^2 + y^2} \le \frac{9}{10}$. The function $v \mapsto v^{-\frac{1}{4}} e^{\frac{\pi}{12} \frac{1}{y} v}$ is bounded in the range $\left[\frac{\sqrt{3}}{2} y, \frac{9}{10}\right]$ by its value at $\frac{9}{10}$, so
\eq
\abs{\frac{f_t(z)}{f_t(iy)}}
< \prn{\frac{9}{10}}^{-\frac{1}{4}} e^{-\frac{\pi}{12} \frac{1}{y} \prn{1 - \frac{9}{10} - \frac{12}{\pi} \frac{1}{265}} + \frac{1}{228}}
< e^{-\frac{1}{70} \frac{1}{y}}
\qe
as claimed.

Now, suppose $t y \le 1$, so that by Lemma \ref{lma:eta-bound-1},
\eq
\eta(ity)^t \ge (ty)^{-\frac{t}{2}} e^{-\frac{\pi}{12} \frac{1}{y} - 1.01 e^{-2\pi} t},
\qe
and hence
\eq
\abs{\frac{f_t(z)}{f_t(iy)}}
\le \abs{\eta(t z)}^t (ty)^{\frac{t}{2}} (y \Im \g z)^{-\frac{1}{4}} e^{\frac{\pi}{12} \Im \g z + \frac{1}{528} t}.
\qe

We will split up the remaining argument into two further cases, depending on whether $\Im \g z \le \frac{\sqrt{3}}{2} t$ or $\Im \g z \ge \frac{\sqrt{3}}{2} t$.

In the first case, when $\Im \g z \le \frac{\sqrt{3}}{2} t$, we may use Lemma \ref{lma:eta-bound-2} to bound $\eta(t z)$ by $\frac{7}{9} (ty)^{-\frac{1}{4}}$ and obtain
\eq
\abs{\frac{f_t(z)}{f_t(iy)}}
\le \prn{\frac{7}{9}}^t (ty)^{\frac{t}{4}} (y \Im \g z)^{-\frac{1}{4}} e^{\frac{\pi}{12} \Im \g z + \frac{1}{528} t},
\qe
which is maximized in the range $\frac{\sqrt{3}}{2} \le \Im \g z \le \frac{\sqrt{3}}{2} t$ when $\Im \g z = \frac{\sqrt{3}}{2} t$, where it is equal to
\eq
\prn{\frac{7}{9}}^{t} (ty)^{\frac{t-1}{4}} \prn{\frac{\sqrt{3}}{2}}^{-\frac{1}{4}} e^{\prn{\frac{\sqrt{3}}{24} \pi + \frac{1}{528}} t}
< 1.04 (ty)^{\frac{t-1}{4}} e^{-\frac{1}{45} t}.
\qe

Let us consider now the second case, where $\Im \g z \ge \frac{\sqrt{3}}{2} t$. Let $\g' = \sMatrix{* & * \\ c/(t,c) & d t / (t,c)} \in \SL_2(\Z)$, which exists since $\frac{c}{(t,c)}$ and $\frac{d t}{(t,c)}$ are coprime. Then,
\eq
\abs{\eta(t z)} = \prn{\frac{ty}{\Im \g' t z}}^{-\frac{1}{4}} \abs{\eta(\g' t z)},
\qe
and
\eq
\Im \g' t z = \frac{t y}{\abs{\frac{c}{(t,c)} t z + \frac{d t}{(t,c)}}^2} = \frac{(t,c)^2}{t} \Im \g z \ge \frac{\sqrt{3}}{2}.
\qe
Therefore,
\eq
\abs{\eta(t z)} \le \sqrt{\frac{(t,c)}{t}} \prn{\frac{y}{\Im \g z}}^{-\frac{1}{4}} e^{-\frac{\pi}{12} \frac{(t,c)^2}{t} \Im \g z + 1.01 e^{-\sqrt{3} \pi}},
\qe
so
\eq \label{eqn:minor-arc-second-case}
\abs{\frac{f_t(z)}{f_t(iy)}}
\le (t,c)^{\frac{t}{2}} (y \Im \g z)^{\frac{t-1}{4}} e^{-\frac{\pi}{12} ((t,c)^2 - 1) \Im \g z + \frac{1}{159} t}.
\qe
Suppose first that $(t,c) > 1$. The function $v \mapsto v^{\frac{t-1}{4}} e^{-\frac{\pi}{12} ((t,c)^2-1) v}$ is maximized at $v = \frac{3}{\pi} \frac{t-1}{(t,c)^2-1} < \frac{\sqrt{3}}{2} t$ and decreasing for larger $v$, so the right-hand side of (\ref{eqn:minor-arc-second-case}) is bounded in the range $\Im \g z \ge \frac{\sqrt{3}}{2} t$ by its value when $\Im \g z = \frac{\sqrt{3}}{2} t$. Hence,
\aleq
\abs{\frac{f_t(z)}{f_t(iy)}}
&\le (t,c)^{\frac{t}{2}} \prn{\frac{\sqrt{3}}{2} t y}^{\frac{t-1}{4}} e^{-\frac{\sqrt{3}}{24} \pi ((t,c)^2 - 1) t + \frac{1}{159} t} \\
&< 2^{\frac{t}{2}} \prn{\frac{\sqrt{3}}{2} t y}^{\frac{t-1}{4}} e^{-\frac{\sqrt{3}}{8} \pi t + \frac{1}{159} t} \\
&< \prn{t y}^{\frac{t-1}{4}} e^{-\frac{1}{3} t}.
\alqe

What remains is the case $(t,c) = 1$, where the right-hand side of (\ref{eqn:minor-arc-second-case}) collapses to become $(y \Im \g z)^{\frac{t-1}{4}} e^{\frac{1}{159} t}$. Letting
\eq
\mathfrak{m}(c, d)
= \set{x \mid \Im \g z \ge \frac{\sqrt{3}}{2} t}
= \set{x \mid \abs{\frac{x + \frac{d}{c}}{y}} < \sqrt{\frac{2}{\sqrt{3}} \frac{1}{c^2 ty} - 1}}
\qe
denote the relevant portion of the interval, we have
\eq
\int_{\frac{y}{3} \le \abs{x} \le \frac{1}{2}} \abs{\frac{f_t(z)}{f_t(iy)}} dx
\le \int_{\substack{\mathfrak{m}(1, 0) \\ \abs{x} > \frac{y}{3}}} \abs{\frac{f_t(z)}{f_t(iy)}} dx + \sum_{c=2}^\infty \sum_{\substack{d=0 \\ (c, d) = 1}}^{c-1} \int_{\mathfrak{m}(c, d)} \abs{\frac{f_t(z)}{f_t(iy)}} dx
+ 1.04 (t y)^{\frac{t-1}{4}} e^{-\frac{1}{45} t}.
\qe

Letting $u = \frac{x + \frac{d}{c}}{y}$, we find
\aleq
\int_{\mathfrak{m}(c, d)} \abs{\frac{f_t(z)}{f_t(iy)}} dx
&\le e^{\frac{1}{159} t} \int_{\mathfrak{m}(c, d)} \prn{\frac{y^2}{\abs{cz+d}^2}}^{\frac{t-1}{4}} dx \\
&= y c^{\frac{1-t}{2}} e^{\frac{1}{159} t} \int_{\abs{u} < \sqrt{\frac{2}{\sqrt{3}} \frac{1}{c^2 ty} - 1}} (u^2 + 1)^{\frac{1-t}{4}} du \\
&< 0.2 y c^{\frac{1-t}{2}} e^{\frac{1}{159} t}, \\
\alqe
so
\aleq
\sum_{c=2}^\infty \sum_{\substack{d=0 \\ (c, d) = 1}}^{c-1} \int_{\mathfrak{m}(c, d)} \abs{\frac{f_t(z)}{f_t(iy)}} dx
< 0.2 y \prn{\z\prn{\tfrac{t-3}{2}} - 1} e^{\frac{1}{159} t}
< y e^{- \frac{1}{3} t}.
\alqe

Letting $v = 1 + (\frac{x}{y})^2$, we also have
\aleq
\int_{\substack{\mathfrak{m}(1, 0) \\ \abs{x} > \frac{y}{3}}} \abs{\frac{f_t(z)}{f_t(iy)}} dx
&\le e^{\frac{1}{159} t} \int_{\substack{\mathfrak{m}(1, 0) \\ \abs{x} > \frac{y}{3}}} \prn{\frac{y^2}{\abs{z}^2}}^{\frac{t-1}{4}} dx \\
&= y e^{\frac{1}{159} t} \int_{\frac{10}{9}}^{\frac{2}{\sqrt{3}} \frac{1}{ty}} v^{\frac{1-t}{4}} \frac{1}{\sqrt{v-1}} dv \\
&< 3 y e^{\frac{1}{159} t} \int_{\frac{10}{9}}^{\infty} v^{\frac{1-t}{4}} dv \\
&< 0.02 y e^{-\frac{1}{50} t}.
\alqe

Altogether, we have
\eq
\int_{\frac{y}{3} \le \abs{x} \le \frac{1}{2}} \abs{\frac{f_t(z)}{f_t(iy)}} dx
< 0.02 y e^{-\frac{1}{50} t}
+ y e^{- \frac{1}{3} t}
+ 1.04 (t y)^{\frac{t-1}{4}} e^{-\frac{1}{45} t}
< y e^{-\frac{1}{70} t},
\qe
as claimed.
\end{proof}



\section{Technical lemmas}
\label{sec:technical-lemmas}
\subsection{Explicit Gaussian integrals} \label{ssec:explicit-integrals}

In \S\ref{sec:asymptotic-formulas}, we will use the following explicit bound for the error term in a Gaussian integral. The proof is elementary, but somewhat long.

\begin{lma} \label{lma:Gaussian-integral-with-error} 
Suppose $\a > 38$, $\b \in \prn{-\frac{2}{25}, \frac{2}{25}}$, and $\varepsilon \in (-1, 1)$. Then,
\[
\begin{aligned}
\int_{-1/3}^{1/3} e(\b x) e^{-\pi \a x^2 \prn{1 + 2 i \varepsilon x + \err 3 x^2}} dx
&= \frac{1}{\sqrt{\a}} \prn{1 + \err \frac{3.45}{\a}} \\
\abs{\int_{-1/3}^{1/3} x e^{-\pi \a x^2 \prn{1 + 2 i \varepsilon x + \err 3 x^2}} dx}
&\le 2 \a^{-\frac{3}{2}}.
\end{aligned}
\]
\end{lma}
\begin{proof}
Let
\eq
I = \int_{-1/3}^{1/3} e(\b x) e^{-\pi \a x^2 \prn{1 + 2 i \varepsilon x + \err 3 x^2}} dx,
\qe
which we decompose as $I = I_1 - \varepsilon I_2 + \err I_3$ where
\aleq
I_1 &= \int_{-1/3}^{1/3} e(\b x) e^{-\pi \a x^2} dx, \\
I_2 &= \int_{-1/3}^{1/3} 2 \pi i \a x^3 e(\b x) e^{-\pi \a x^2} dx, \quad \text{and} \\
I_3 &= \int_{-1/3}^{1/3} \abs{(1 - 2 \pi i \varepsilon \a x^3) e(\b x) e^{-\pi \a x^2} - e(\b x) e^{-\pi \a x^2 \prn{1 + 2 i \varepsilon x + \err 3 x^2}}} dx.
\alqe

The main term $I_1$ is simply
\aleq
I_1
&= \int_{-\infty}^\infty e(\b x) e^{-\pi \a x^2} dx + \err 2 \int_{1/3}^\infty e^{-\pi \a x^2} dx \\
&= e^{-\pi \frac{\b^2}{\a}} \frac{1}{\sqrt{\a}} + \err 2 \int_{1/3}^\infty \frac{2\pi \a x}{2 \pi \a \frac{1}{3}} e^{-\pi \a x^2} dx \\
&= e^{-\pi \frac{\b^2}{\a}} \frac{1}{\sqrt{\a}} + \err \frac{3}{\pi \a} e^{-\frac{\pi}{9} \a}.
\alqe

Similarly, the term $I_2$ is
\aleq
I_2
&= 2\pi \a \int_{-\infty}^\infty i x^3 e(\b x) e^{-\pi \a x^2} dx + \err 4 \pi \a \int_{1/3}^\infty x^3 e^{-\pi \a x^2} dx \\
&= e^{-\pi \frac{\b^2}{\a}} \frac{1}{\sqrt{\a}} \prn{-3 \frac{\b}{\a} + 2 \pi \frac{\b^3}{\a^2}} + \err \prn{\frac{2}{9} + \frac{2}{\pi \a}} e^{-\frac{\pi}{9} \a}
\alqe

As for $I_3$, using Taylor's theorem in the form $e^z = 1 + z e^{s z} = 1 + z + \frac{1}{2} z^2 e^{s' z}$ for some $s, s' \in (0, 1)$, we find that the integrand in $I_3$ is equal to
\aleq
&e^{-\pi \a x^2} \abs{(1 - 2 \pi i \varepsilon \a x^3) - e^{-\pi \a x^2 \prn{2 i \varepsilon x + \err 3 x^2}}} \\
&\le  e^{-\pi \a x^2} \prn{\abs{(1 - 2 \pi i \varepsilon \a x^3) - e^{-2 \pi i \varepsilon \a x^3}} + \abs{e^{-2 \pi i \varepsilon \a x^3} - e^{-\pi \a x^2 \prn{2 i \varepsilon x + \err 3 x^2}}}} \\
&\le e^{-\pi \a x^2} \prn{2 \pi^2 \a^2 x^6 + 3\pi \a x^4 e^{3 \pi \a x^4}}.
\alqe
Therefore,
\aleq
I_3
&\le \int_{-1/3}^{1/3} \prn{2 \pi^2 \a^2 x^6 e^{-\pi \a x^2} + 3\pi \a x^4 e^{-\pi \a x^2 (1 - 3 x^2)}} dx \\
&\le 2 \pi^2 \a^2 \int_{-\infty}^\infty x^6 e^{-\pi \a x^2} dx + 3\pi \a \int_{-\infty}^\infty x^4 e^{-\frac{2}{3} \pi \a x^2} dx \\
&= \prn{\frac{15}{4 \pi} + \frac{(3/2)^{\frac{9}{2}}}{\pi}} \a^{-\frac{3}{2}}.
\alqe
Altogether, we have
\aleq
I
&= \frac{e^{-\pi \frac{\b^2}{\a}}}{\sqrt{\a}} \prn{1 + 3 \varepsilon \frac{\b}{\a} - 2 \pi \varepsilon \frac{\b^3}{\a^2}}
+ \err e^{-\frac{\pi}{9} \a} \prn{\frac{5}{\pi \a} + \frac{2}{9}} + \err \prn{\frac{15}{4 \pi} + \frac{(3/2)^{\frac{9}{2}}}{\pi}} \a^{-\frac{3}{2}} \\
&= \frac{1}{\sqrt{\a}} \prn{1 + \err \frac{3.45}{\a}},
\alqe
as claimed.

We will now treat the second integral in the lemma. Let
\eq
J = \int_{-1/3}^{1/3} x e^{-\pi \a x^2 \prn{1 + 2 i \varepsilon x + \err 3 x^2}} dx.
\qe
Since $x e^{-\pi \a x^2}$ is odd, we can write
\eq
J = \int_{-1/3}^{1/3} x \prn{e^{-\pi \a x^2 \prn{1 + 2 i \varepsilon x + \err 3 x^2}} - e^{-\pi \a x^2}} dx.
\qe
The integrand has magnitude bounded by
\aleq
&\abs{x} e^{-\pi \a x^2} \abs{e^{-2\pi i \varepsilon \a x^3 + \err 3 \pi \a x^4} - 1} \\
&\le \abs{x} e^{-\pi \a x^2} \abs{-2 \pi i \varepsilon \a x^3 + \err 3 \pi \a x^4} e^{3 \pi \a x^4} \\
&\le \prn{2 \pi \a x^4 + 3 \pi \a \abs{x}^5} e^{-\pi \a x^2 (1 - 3 x^2)},
\alqe
so we find
\aleq
\abs{J}
&\le \int_{-1/3}^{1/3} \prn{2 \pi \a x^4 + 3 \pi \a \abs{x}^5} e^{-\pi \a x^2 (1 - 3 x^2)} dx \\
&\le \int_{-\infty}^\infty \prn{2 \pi \a x^4 + 3 \pi \a \abs{x}^5} e^{-\frac{2}{3} \pi \a x^2} dx \\
&= \frac{(3/2)^{\frac{7}{2}}}{\pi} \a^{-\frac{3}{2}} + \frac{81}{4 \pi^2} \a^{-2} \\
&< 2 \a^{-\frac{3}{2}}.
\alqe
\end{proof}

\subsection{The functions \texorpdfstring{$\func_k$}{funck}} \label{ssec:func-bounds}

Recall the functions
\eq
\func_k(z) = - \frac{z^{k+1}}{2 \pi i} \prn{\dv{}{z}}^k \log \eta(z)
\qe
from the introduction (\ref{eqn:func-def}), which appear in the Taylor expansion
\eq
\frac{1}{2\pi i} \log f_t(z)
= \sum_{k=0}^\infty \frac{x^k}{k!} \frac{\func_k(iy) - \func_k(ity)}{(i y)^{k+1}}.
\qe
Our main goal in this section is to bound the third and fourth terms in this Taylor expansion in terms of the second.

Using the functional equation $\eta(z) = (-iz)^{-\frac{1}{2}} \eta(-1/z)$, we find two different expansions for $\func_k(z)$, useful depending on whether $\Im z$ or $\Im \frac{-1}{z}$ is large. We have
\aleq \label{eqn:func-fourier}
\func_k(z)
&= \sum_{n=1}^\infty z^{k+1} (2\pi i n)^{k-1} \s(n) e(nz)
- \begin{cases}
    z^2/24 & \text{if $k = 0,1$} \\
    0 & \text{if $k \ge 2$}
\end{cases} \\
&= 
\sum_{n=1}^\infty P_k\prn{\frac{2\pi i n}{z}} \s(n) e(\tfrac{-n}{z})
+ \frac{(-1)^k k!}{24}
+ \frac{z}{4\pi i} \begin{cases}
    \log(-iz) & \text{if $k = 0$} \\
    (-1)^{k-1} (k-1)! & \text{if $k \ge 1$}
\end{cases}
\alqe
where $P_k(r)$ is the degree $k-1$ polynomial (Laurent in the case $k = 0$) defined by the recurrence $P_0(r) = r^{-1}$ and
\eq
P_k(r) = (r - k) P_{k-1}(r) - r P_{k-1}'(r).
\qe

We will also need expansions for $\func_k'(z)$. We have
\aleq \label{eqn:func-prime-fourier}
\func_k'(z)
&= \sum_{n=1}^\infty \prn{(2 \pi i n z)^{k+1} + (k+1) (2 \pi i n z)^k} \frac{\s(n)}{2\pi i n} e(n z)
- \begin{cases}
    z/12 & \text{if $k = 0,1$} \\
    0 & \text{if $k \ge 2$}
\end{cases} \\
&= \sum_{n=1}^\infty Q_k\prn{\frac{2\pi i n}{z}} \frac{\s(n)}{2\pi i n} e(\tfrac{-n}{z})
+ \frac{1}{4\pi i} \begin{cases}
    1 + \log(-iz) & \text{if $k = 0$} \\
    (-1)^{k-1} (k-1)! & \text{if $k \ge 1$}
\end{cases}
\alqe
where $Q_k(r)$ is the degree $k+1$ polynomial defined by the recurrence $Q_0(r) = r+1$, and
\aleq
Q_k(r)
&= (r-k+1) Q_{k-1}(r) - r Q_{k-1}'(r) \\
&= r^2 (P_k(r) - P_k'(r)).
\alqe

The first few values of $P_k(r)$ and $Q_k(r)$ are given in the table (\ref{tbl:polynomials}) below.
\eq
\begin{array}{c|c|c} \label{tbl:polynomials}
k & P_k(r) & Q_k(r) \\ \hline
0 & r^{-1} & r + 1 \\
1 & 1 & r^2 \\
2 & r-2 & r^3 - 3r^2 \\
3 & r^2 - 6r + 6 & r^4 - 8r^3 + 12 r^2 \\
4 & r^3 - 12 r^2 + 36 r - 24 & r^5 - 15 r^4 + 60 r^3 - 60 r^2
\end{array}
\qe

As a first observation, we can see directly from (\ref{eqn:func-prime-fourier}) and (\ref{tbl:polynomials}) that
\eq \label{eqn:func_2-functional-equation}
\func_2'(z) + \func_2'(\tfrac{-1}{z}) = -\frac{1}{4\pi i}.
\qe

\begin{lma} \label{lma:func_2-bound}
\leavevmode
\begin{enumerate}[(i)]
    \item \label{item:func_2-prime-increasing} The function $y \mapsto i \func_2'(iy)$ is a monotonically increasing bijection $(0, \infty) \to (-\frac{1}{4\pi}, 0)$.
    
    \item \label{item:func_2-decreasing} The function $y \mapsto \func_2(iy)$ is a monotonically decreasing bijection $(0, \infty) \to (0, \frac{1}{12})$.
\end{enumerate}
Suppose $0 < y \le \frac{1}{10}$ and $t > 1$.
\begin{enumerate}[(i)]
\setcounter{enumi}{2}
    \item \label{item:func_2-diff-small} If $t y \le 1$, then
    \[
    \frac{1}{8 \pi} < \frac{\func_2(iy) - \func_2(ity)}{ty-y} < \frac{1}{4 \pi}.
    \]

    \item \label{item:func_2-diff-big} If $t y \ge 1$, then
    \[
    \frac{1}{16} < \func_2(iy) - \func_2(ity) < \frac{1}{12}.
    \]
\end{enumerate}
\end{lma}

\begin{proof}
\leavevmode
\begin{enumerate}[(i)]
    \item From (\ref{eqn:func-prime-fourier}), we see  $\lim_{y\to0} i \func_2'(iy) = -\frac{1}{4\pi}$ and $\lim_{y\to\infty} i \func_2'(iy) = 0$, so it remains to show that $i \func_2'(iy)$ is monotonically increasing. By (\ref{eqn:func_2-functional-equation}), we may assume $y \ge 1$, in which case each term in the sum
    \eq
    i \func_2'(iy)
    = \sum_{n=1}^\infty \prn{-(2 \pi n y)^3 + 3 (2\pi n y)^2} \frac{\s(n)}{2 \pi n} e^{-2\pi n y}
    \qe
    is monotonically increasing.
    
    \item From (\ref{eqn:func-fourier}), we see $\lim_{y\to0} \func_2(iy) = \frac{1}{12}$ and $\lim_{y\to\infty} \func_2(iy) = 0$, and $\func_2(iy)$ is monotonically decreasing by part (i).

    \item By the mean value theorem, $\frac{\func_2(iy) - \func_2(iv)}{y - ty} = i \func_2'(i v)$ for some $v \in (y, ty) \subset (0, 1]$. By part (i),
    \eq
    -\frac{1}{4\pi} < i \func_2'(i v) < i \func_2'(i) = -\frac{1}{8 \pi}.
    \qe

    \item By part (ii),
    \eq
    0.0635 \approx \func_2\prn{\frac{i}{10}} - \func_2(i) \le \func_2(iy) - \func_2(ity) < \frac{1}{12} - 0.
    \qe
    
\end{enumerate}
\end{proof}

\begin{lma} \label{lma:func_3-bound}
\leavevmode
\begin{enumerate}[(i)]
    \item For all $y > 0$, $0 < i \func_3'(iy) < \frac{3}{4\pi}$.

    \item The function $y \mapsto \func_3(iy)$ is a monotonically increasing bijection $(0, \infty) \to (-\frac{1}{4}, 0)$.
    
    \item If $0 < y \le \frac{1}{10}$ and $t > 1$, then
    \[
    \abs{\frac{\func_3(iy) - \func_3(ity)}{\func_2(iy) - \func_2(ity)}} < 6.
    \]
\end{enumerate}
\end{lma}

\begin{proof}
\leavevmode
\begin{enumerate}[(i)]

    \item First, suppose $y \le 1$. The function $Q_3(r) e^{-r} = r^2 (r-2) (r-6) e^{-r}$ is positive and bounded by $0.3$ for $r > 6$, and decreasing for $r > 8$, so
    \aleq
    0 < i \func_3'(iy)
    &= \frac{1}{2\pi} + \sum_{n=1}^\infty Q_3\prn{\frac{2\pi n}{y}} \frac{\s(n)}{2\pi n} e^{-\frac{2\pi n}{y}} \\
    &< \frac{1}{2\pi} + \frac{0.3}{2\pi} + \sum_{n=2}^\infty Q_3(2 \pi n) \frac{n}{2\pi} e^{-2\pi n} \\
    &< \frac{3}{4\pi}.
    \alqe
    Now, suppose $y \ge 1$. The function $(r^4 - 4 r^3) e^{-r}$ is positive and decreasing for $r > 6$, so
    \aleq
    0 < i \func_3'(iy)
    &= \sum_{n=1}^\infty \prn{(2\pi n y)^4 - 4 (2 \pi n y)^3} \frac{\s(n)}{2\pi n} e^{-2\pi n y} \\
    &< \sum_{n=1}^\infty \prn{(2\pi n)^4 - 4 (2\pi n)^3} \frac{n}{2\pi} e^{-2\pi n} \\
    &< \frac{3}{4\pi}.
    \alqe

    \item From (\ref{eqn:func-fourier}), we see $\lim_{y\to0} \func_3(iy) = -\frac{1}{4}$ and $\lim_{y\to\infty} \func_3(iy) = 0$, and $\func_3(iy)$ is monotonically increasing by part (i).

    \item If $ty \le 1$, then by part \ref{item:func_2-diff-small} of Lemma \ref{lma:func_2-bound},
    \eq
    \abs{\frac{\func_3(iy) - \func_3(ity)}{\func_2(iy) - \func_2(ity)}}
    < 8 \pi \abs{\frac{\func_3(iy) - \func_3(ity)}{y - ty}}
    = 8 \pi \abs{\func_3'(isy)}
    < 6.
    \qe
    for some $s \in (1, t)$.

    If $ty \ge 1$, then by part \ref{item:func_2-diff-big} of Lemma \ref{lma:func_2-bound},
    \eq
    \abs{\frac{\func_3(iy) - \func_3(ity)}{\func_2(iy) - \func_2(ity)}}
    < \frac{1/4 - 0}{1/16} = 4.
    \qe
    
\end{enumerate}

\end{proof}

\begin{lma} \label{lma:func_4-bound}
Suppose $z = x + iy$ with $0 < y \le \frac{1}{10}$ and $\abs{x} < \frac{y}{3}$. For any $t > 1$,
\[
\abs{\frac{\func_4(z) - \func_4(tz)}{\func_2(iy) - \func_2(ity)}} < 36.
\]
\end{lma}

\begin{proof}
We will treat separately the cases $t y \le \frac{9}{10}$ and $t y \ge \frac{9}{10}$.

First, suppose $t y \le \frac{9}{10}$. By the mean value theorem,
\eq
\abs{\func_4(z) - \func_4(tz)} = (t-1) \abs{z} \abs{\func_4'(sz)}
\qe
for some $s \in (1, t)$. Note that
\eq
\func_4'(sz) = -\frac{3}{2\pi i} + \sum_{n=1}^\infty Q_4\prn{\frac{2\pi i n}{sz}} \frac{\s(n)}{2 \pi i n} e^{-\frac{2 \pi i n}{sz}}.
\qe
Let $r = \frac{2\pi i}{sz}$, and note that $\abs{\Im r} < \frac{\Re r}{3}$ and $\Re r > 2 \pi$. When $u$ is fixed, the function $\abs{Q_4(u+iv)}$ is increasing in $\abs{v}$ since $Q_4(r)$ has real roots. Moreover, the function $\abs{Q_4(u(1+\frac{i}{3}))} e^{-u}$ is decreasing in $u$ for $u > 3$. Therefore, in the region $\abs{\Im r} < \frac{\Re r}{3}$ and $\Re r > 2 \pi$, the function $\abs{Q_4(n r) e^{-n r}}$ is bounded by $\abs{Q_4(2 \pi n(1+\frac{i}{3}))} e^{-2 \pi n}$, which means
\eq
\abs{\func_4'(sz)}
\le \frac{3}{2\pi} + \sum_{n=1}^\infty \abs{Q_4(2 \pi n (1 + i/3)} \frac{\s(n)}{2\pi n} e^{-2\pi n}
< \frac{4}{\pi}.
\qe
Therefore, by part \ref{item:func_2-diff-small} of Lemma \ref{lma:func_2-bound},
\eq
\abs{\frac{\func_4(z) - \func_4(tz)}{\func_2(iy) - \func_2(ity)}}
< \frac{4 / \pi}{1/(8\pi)} \frac{|z|}{y}
< 32 \sqrt{\frac{10}{9}}
< 36.
\qe

Now, suppose $t y \ge \frac{9}{10}$. Then,
\aleq
\abs{\func_4(tz)}
&= \abs{\sum_{n=1}^\infty z^5 (2 \pi i n)^3 \s(n) e(nz)}
&\le \prn{\frac{10}{9}}^{\frac{5}{2}} \sum_{n=1}^\infty \prn{\frac{9}{10}}^5 (2 \pi n)^3 n^2 e^{-2\pi \frac{9}{10} n}
&< 1,
\alqe
and since $\abs{P_4\prn{\frac{2\pi i n}{z}}} < 10^5 (\frac{n}{y})^3$,
\aleq
\abs{\func_4(z)}
&\le \abs{1 - \frac{3z}{2\pi i}} + \sum_{n=1}^\infty \abs{P_4\prn{\frac{2\pi i n}{z}} \s(n) e(\tfrac{-n}{z})}
\le 1 + 10^5 \sum_{n=1}^\infty n^5 e^{-18\pi n}
&\le 1.01.
\alqe
Therefore, by part \ref{item:func_2-decreasing} of Lemma \ref{lma:func_2-bound},
\eq
\abs{\frac{\func_4(z) - \func_4(tz)}{\func_2(iy) - \func_2(ity)}}
< \frac{2.01}{\func_2(\frac{1}{10} i) - \func_2(\frac{9}{10} i)}
< 36.
\qe
\end{proof}

\section{Asymptotic formulas for \texorpdfstring{$c_t(N)$}{ct(N)}}
\label{sec:asymptotic-formulas}
We are now ready to carry out the argument sketched in \S \ref{sec:sketch}. Recall that $M = N + \frac{t^2-1}{24}$.

\subsection{Proof of part \ref{item:y} of Theorem \ref{thm:t-core-asymptotic}} \label{ssec:y-def}

By (\ref{eqn:func-fourier}), we have
\eq
\func_1(iy)
= \frac{y^2}{24} - \sum_{n=1}^\infty y^2 \s(n) e^{-2\pi n y}
= -\frac{1}{24} + \frac{y}{4\pi} + \sum_{n=1}^\infty \s(n) e^{-2\pi n / y},
\qe
which means
\eq
\lim_{y\to\infty} \frac{\func_1(ity) - \func_1(iy)}{y^2} = \frac{t^2-1}{24} \quad \text{and} \quad
\lim_{y\to0^+}    \frac{\func_1(ity) - \func_1(iy)}{y^2} = \infty,
\qe
and hence there exists $y > 0$ such that $\frac{\func_1(ity) - \func_1(iy)}{y^2} = M$.

Applying part \ref{item:func_2-decreasing} of Lemma \ref{lma:func_2-bound}, we have
\eq \label{eqn:y-unique}
\dv{}{y} \prn{\frac{\func_1(ity) - \func_1(iy)}{y^2}} = \frac{\func_2(ity) - \func_2(iy)}{y^3} < 0,
\qe
so this value of $y$ is unique.

Moreover,
\eq
-\frac{1}{24} + \frac{y}{4\pi} < \func_1(iy) < \frac{y^2}{24},
\qe
so
\eq \label{eqn:func_1-range-1}
\frac{\func_1(ity) - \func_1(iy)}{y^2} < \frac{t^2}{24} + \frac{1}{24 y^2} - \frac{1}{4\pi y}.
\qe

By (\ref{eqn:func-prime-fourier}), we also have
\eq
i \func_1'(iy) = \frac{1}{4\pi} + \sum_{n=1}^\infty \prn{\frac{2\pi n}{y}}^2 \frac{\s(n)}{2\pi n} e^{-2\pi n / y} > \frac{1}{4\pi},
\qe
so by the mean value theorem,
\eq \label{eqn:func_1-range-2}
\frac{\func_1(ity) - \func_1(iy)}{y^2} > \frac{t-1}{4\pi y}.
\qe

Together, (\ref{eqn:func_1-range-1}) and (\ref{eqn:func_1-range-2}) imply that if $\frac{\func_1(ity) - \func_1(iy)}{y^2} = M$, then
\eq
\frac{t-1}{4\pi M} < y < \frac{1}{\frac{3}{\pi} + \sqrt{24 N - 1 + \frac{9}{\pi^2}}}.
\qe

\subsection{Proof of part \ref{item:asymp} of Theorem \ref{thm:t-core-asymptotic}} \label{ssec:asymptotic-proof}

Since $c_t(N) = p(N)$ if $t > N$, we may assume $N$ is sufficiently large. By Lulov and Pittel's formula (\ref{eqn:anderson}), we may also assume $t$ is sufficiently large. Moreover, parts \ref{item:func_2-diff-small} and \ref{item:func_2-diff-big} of Lemma \ref{lma:func_2-bound} and (\ref{eqn:func_1-range-1}) imply that if $\abs{\frac{\func_1(ity) - \func_1(iy)}{y^2} - M} \ll \frac{1}{y}$, then
\eq
\frac{y}{\func_2(iy) - \func_2(ity)} \ll \frac{1}{\min(t, \frac{1}{y})} \ll \frac{1}{\min(t, \sqrt{N})}.
\qe
Therefore, the following theorem implies (and makes explicit) part \ref{item:asymp} of Theorem \ref{thm:t-core-asymptotic}.

\begin{thm} \label{thm:t-core-asymptotic-explicit}
Choose $y$ so that
\[
\abs{\frac{\func_1(ity) - \func_1(iy)}{y^2} - M} < \frac{2}{25 y}
\]
(e.g. choosing $y$ as in \S\ref{ssec:y-def} sets this quantity to $0$), and suppose $\min(t, \frac{1}{y}) \ge 1000$. Then,
\[
c_t(N) = \frac{y^{\frac{3}{2}} e^{2\pi M y} f_t(iy)}{\sqrt{\func_2(iy) - \func_2(ity)}}
\prn{
1 + \err \frac{3.5 y}{\func_2(iy) - \func_2(ity)}
}.
\]
\end{thm}

\begin{proof}

Consider the integral
\aleq \label{eqn:t-core-integral}
c_t(N)
&= \int_{-1/2}^{1/2} e(-M z) f_t(z) dx \\
&= e^{2\pi M y} f_t(iy) \prn{\int_{-y/3}^{y/3} e\prn{-M x + \frac{1}{2\pi i} \log \frac{f_t(z)}{f_t(iy)}} dx + \err y e^{-\frac{1}{70} \min(t, \frac{1}{y})}},
\alqe
where we have applied Proposition \ref{prop:minor-arc-bound} to obtain the second line from the first.

We have the Taylor expansion
\aleq \label{eqn:taylor-expansion}
\frac{1}{2\pi i} \log \frac{f_t(z)}{f_t(iy)}
&= x \frac{\func_1(iy) - \func_1(ity)}{(iy)^2}
+ \frac{x^2}{2} \frac{\func_2(iy) - \func_2(ity)}{(iy)^3} \\
&+ \frac{x^3}{6} \frac{\func_3(iy) - \func_3(ity)}{(iy)^4}
+ \frac{x^4}{24} \frac{\func_4(z') - \func_4(tz')}{z'^5}
\alqe
for some $z' = x' + iy$ with $x'$ between $0$ and $x$. Lemmas \ref{lma:func_3-bound} and \ref{lma:func_4-bound} allow us to bound the third and fourth terms above in terms of the quadratic term, and we can rewrite (\ref{eqn:taylor-expansion}) as
\eq
\frac{1}{2\pi i} \log \frac{f_t(z)}{f_t(iy)}
= x \frac{\func_1(iy) - \func_1(ity)}{(iy)^2}
+ \frac{x^2}{2} \frac{\func_2(iy) - \func_2(ity)}{(iy)^3}
\prn{1 + 2 i \varepsilon \frac{x}{y} + \err 3 \frac{x^2}{y^2}}
\qe
where $\varepsilon \in (-1, 1)$.

Therefore, if we let
\eq
\a = \frac{\func_2(iy) - \func_2(ity)}{y} \quad \text{and} \quad
\b = y \prn{\frac{\func_1(iy) - \func_1(ity)}{(iy)^2} - M},
\qe
which satisfy
\eq
\frac{1}{26} \le \frac{\a}{\min(t, 1/y)} \le \frac{1}{12} \quad \text{and} \quad
\abs{\b} < \frac{2}{25}
\qe
by parts \ref{item:func_2-diff-small} and \ref{item:func_2-diff-big} of Lemma \ref{lma:func_2-bound} and the assumption in Theorem \ref{thm:t-core-asymptotic-explicit}, then
\eq \label{eqn:reduction-to-gaussian}
e\prn{-Mx + \frac{1}{2\pi i} \log \frac{f_t(z)}{f_t(iy)}}
= e\prn{\b x/y} e^{-\a (x/y)^2 \prn{1 + 2 i \e (x/y) + \err 3 (x/y)^2}},
\qe
and by Lemma \ref{lma:Gaussian-integral-with-error}, making the change of variables $u = x/y$, the integral in the second line of (\ref{eqn:t-core-integral}) is simply
\eq
y \int_{-1/3}^{1/3} e(\b u) e^{- \pi \a u^2 \prn{1 + 2 i \varepsilon u + \err 3 u^2}} du
= \frac{y}{\sqrt{\a}} \prn{1 + \frac{3.45}{\a}}.
\qe

Together with (\ref{eqn:t-core-integral}), this shows
\aleq
c_t(N)
&= e^{2\pi M y} f_t(iy) \prn{\frac{y}{\sqrt{\a}} \prn{1 + \err \frac{3.45}{\a}} + \err y e^{-\frac{1}{70} \min(t, \frac{1}{y})}} \\
&= \frac{y^{\frac{3}{2}} e^{2\pi M y} f_t(iy)}{\sqrt{\func_2(iy) - \func_2(ity)}} \prn{1 + \err \frac{3.45}{\a} + \err \sqrt{\a} e^{-\frac{1}{70} \min(t, \frac{1}{y})}} \\
&= \frac{y^{\frac{3}{2}} e^{2\pi M y} f_t(iy)}{\sqrt{\func_2(iy) - \func_2(ity)}} \prn{1 + \err \frac{3.5}{\a}},
\alqe
which completes the proof of Theorem \ref{thm:t-core-asymptotic-explicit}.
\end{proof}

Note also that
\eq \label{eqn:f_t-abs-integral}
\int_{-y/3}^{y/3} \abs{f_t(z)} dx
\le f_t(iy) \int_{-y/3}^{y/3} e^{-\pi \a (x/y)^2 \prn{1 + \err 3(x/y)^2}} dx
\le \sqrt{\frac{3}{2\a}} y f_t(iy),
\qe
a fact which will be useful in \S\ref{ssec:diff-proof}.

\subsection{Proof of Corollary \ref{cly:t-core-asymptotic-middle-t}} \label{ssec:middle-t-proof}

First, note that
\aleq
\lim_{v\to\infty} \prn{\frac{1}{24 v^2} - \frac{1}{24} + \frac{1}{v^2} \func_1(iv)}
&= 0 \quad \text{and} \\
\lim_{v\to0^+}    \prn{\frac{1}{24 v^2} - \frac{1}{24} + \frac{1}{v^2} \func_1(iv)}
&= \infty,
\alqe
and by part \ref{item:func_2-decreasing} of Lemma \ref{lma:func_2-bound},
\eq
\dv{}{v} \prn{\frac{1}{24 v^2} - \frac{1}{24} + \frac{1}{v^2} \func_1(iv)}
= \frac{1}{v^3} \prn{\func_2(iv) - \frac{1}{12}}
< 0.
\qe
Hence, the function $\frac{1}{24 v^2} - \frac{1}{24} + \frac{1}{v^2} \func_1(iv)$ is a monotonically decreasing bijection $(0, \infty) \to (0, \infty)$, and $v$ exists and is unique.

Let $y = \frac{v}{t}$. We begin by showing that the hypotheses of Theorem \ref{thm:t-core-asymptotic-explicit} are satisfied for $t$ sufficiently large. Indeed,
\aleq
\abs{\frac{\func_1(ity) - \func_1(iy)}{y^2} - M}
&= \abs{\frac{1}{y^2} \func_1(iv) + \frac{1}{24 y^2} - \frac{1}{4\pi y} + O\prn{y^{-2} e^{-2\pi/y}} - \frac{t^2}{\kappa} - \frac{t^2-1}{24} + o(t)} \\
&= \abs{\frac{v^2}{y^2} \prn{\frac{1}{24 v^2} - \frac{1}{24} + \frac{1}{v^2} \func_1(iv) - \frac{1}{\kappa}} - \frac{1}{4\pi y}} + o(t) \\
&= \frac{1}{4\pi y} \prn{1 + o(1)}.
\alqe
Therefore, by Theorem \ref{thm:t-core-asymptotic-explicit},
\eq \label{eqn:middle-t-c_t-expression}
c_t(N) = \frac{y}{\sqrt{\a}} e^{2\pi M y} f_t(iy) \prn{1 + O\prn{\frac{1}{\a}}}
\qe
where
\aleq
\a
= \frac{\func_2(iy) - \func_2(ity)}{y}
= \frac{1}{y} \prn{\frac{1}{12} - \func_2(iv)} + O(1).
\alqe

We also have
\aleq
\log f_t(iy)
= \frac{2\pi}{y} \prn{\func_0(iy) - \func_0(ity)}
= \frac{2\pi}{y} \prn{\frac{1}{24} - \func_0(iv)} + \frac{1}{2} \log y + O\prn{y^{-2} e^{-2\pi/y}},
\alqe
so (\ref{eqn:middle-t-c_t-expression}) becomes
\aleq
c_t(N)
&= \frac{y^2}{\sqrt{\frac{1}{12} - \func_2(iv)}} \exp\prn{\frac{2\pi}{y} \prn{\prn{N + \frac{t^2-1}{24}} y^2 + \frac{1}{24} - \func_0(iv)}} \prn{1 + O\prn{y}}.
\alqe

We conclude by noting
\eq
\frac{y^2}{\sqrt{\frac{1}{12} - \func_2(iv)}}
= \frac{\kappa}{t^2 B(\kappa)}
= \frac{1}{B(\kappa) N} \prn{1 + o\prn{y}}
\qe
and
\aleq
\frac{2\pi}{y} \prn{\prn{N + \frac{t^2-1}{24}} y^2 + \frac{1}{24} - \func_0(iv)}
&= 2 \pi \sqrt{\frac{\kappa N}{v^2}} \prn{\frac{v^2}{\kappa} + \frac{v^2}{24} + \frac{1}{24} - \func_0(iv)} + o(1) \\
&= 2\pi \sqrt{A(\kappa) N} + o(1). \\
\alqe

\subsection{Proof of part \ref{item:diff-asymp} of Theorem \ref{thm:t-core-asymptotic}} \label{ssec:diff-proof}

We will prove the following explicit version of part \ref{item:diff-asymp} of Theorem \ref{thm:t-core-asymptotic}, from which part \ref{item:diff-asymp} of Theorem \ref{thm:t-core-asymptotic} immediately follows.

\begin{thm} \label{thm:t-core-diff-asymptotic-explicit}
Choose $y$ as in \S\ref{ssec:y-def}, and suppose $\min(t, \frac{1}{y}) \ge 1000$, and $t y \ge \frac{1}{2}$. Then,
\[
c_{t+1}(N+t) - c_t(N+t)
= c_t(N) \prn{
2\pi t y - 1 + \err t y \prn{705 y + 120 t y e^{-2\pi t y}}}.
\]
\end{thm}

We begin with the following lemma, which will allow us to interpolate between $c_{t+1}(N)$ and $c_t(N)$.

\begin{lma} \label{lma:func_0-bound}
If $ty \ge \frac{1}{2}$, then
\[
\begin{aligned}
\abs{\log \prn{e\prn{-\tfrac{2t+1}{24} z} \frac{f_{t+1}(z)}{f_t(z)}}}
&\le 7.5 \abs{t z} e^{-2\pi t y} \quad \text{and} \\
\abs{\log \prn{e\prn{-\tfrac{2t+1}{24} z} \frac{f_{t+1}(z)}{f_t(z)}} + e(tz) (2 \pi i t z + 1)}
&\le \prn{40 \abs{z} + 22 e^{-2\pi t y}} \abs{t z} e^{-2\pi t y}.
\end{aligned}
\]
\end{lma}
\begin{proof}
By the mean value theorem,
\aleq \label{eqn:func_0-mvt}
\frac{1}{2\pi i} \log \prn{e\prn{-\tfrac{2t+1}{24} z} \frac{f_{t+1}(z)}{f_t(z)}}
&= \frac{1}{z} \prn{\func_0(tz) + \frac{(tz)^2}{24} - \prn{\func_0((t+1)z) + \frac{((t+1) z)^2}{24}}} \\
&= -\func_0'((t+\d) z) - \frac{(t+\d) z}{12}
\alqe
for some $\d \in (0,1)$. Let $w = (t+\d) z$, $u = (t+\d) x$, and $v = (t+\d) y$ (so $w = u + i v$).

Consider the function
\eq
\func_0'(w) + \frac{w}{12}
= \sum_{n=1}^\infty (2 \pi i n w + 1) \frac{\s(n)}{2\pi i n} e(n w)
\qe
and note that if $v \ge \frac{1}{2}$, then
\eq
\abs{\func_0'(w) + \frac{w}{12}}
\le \sum_{n=1}^\infty \abs{2\pi i n w + 1} \frac{\s(n)}{2\pi n} e^{-2\pi n v}
\le \abs{w} \sum_{n=1}^\infty \s(n) e^{-2\pi n v}
< \frac{15}{4 \pi} \abs{w} e^{-2\pi v}.
\qe

Therefore,
\eq
\abs{\log \prn{e\prn{-\tfrac{2t+1}{24} z} \frac{f_{t+1}(z)}{f_t(z)}}}
\le 2 \pi \abs{\func_0'(w) + \frac{w}{12}}
\le 7.5 \abs{w} e^{-2\pi v}
\le 7.5 \abs{t z} e^{-2\pi t y},
\qe
which is the first inequality claimed. For the second, we have
\aleq \label{eqn:func_0-mvt-1}
&{\abs{\log \prn{e\prn{-\tfrac{2t+1}{24} z} \frac{f_{t+1}(z)}{f_t(z)}} + (2 \pi i w + 1) e(w)}} \\
&= 2 \pi \abs{-\func_0^\prime(w) - \frac{w}{12} + (2\pi i w + 1) \frac{1}{2\pi i} e(w)} \\
&= 2 \pi \abs{\sum_{n=2}^\infty (2\pi i n w + 1) \frac{\s(n)}{2\pi i n} e(nw)} \\
&\le 2 \pi \abs{w} \sum_{n=2}^\infty \s(n) e^{-2\pi n v} \\
&< 22 \abs{w} e^{-4 \pi v} \\
&\le 22 \abs{t z} e^{-4 \pi t y}
\alqe
and, applying the mean value theorem,
\eq \label{eqn:func_0-mvt-2}
\abs{(2 \pi i w + 1) e(w) - (2 \pi i t z + 1) e(t z)}
\le 40 \abs{t z^2} e^{-2\pi t y}.
\qe
Adding (\ref{eqn:func_0-mvt-1}) and (\ref{eqn:func_0-mvt-2}) completes the proof.
\end{proof}

\begin{proof}[Proof of Theorem \ref{thm:t-core-diff-asymptotic-explicit}]
Recall that we assume $\min(t, \frac{1}{y}) \ge 1000$ and $ty \ge \frac{1}{2}$.

Consider the integral
\eq
c_{t+1}(N+t) - c_t(N+t)
= \int_{-1/2}^{1/2} e(-(M+t) z) f_t(z) \prn{e\prn{-\tfrac{2t+1}{24} z}\frac{f_{t+1}(z)}{f_t(z)} - 1} dx.
\qe

We will begin by bounding the contribution to this integral where $\abs{x} \ge \frac{y}{3}$. By splitting up into two cases corresponding to $ty \le \frac{1}{2\pi} \log(4t)$ and $t y \ge \frac{1}{2\pi} \log(4t)$, we will show
\eq \label{eqn:diff-away-from-0}
\int_{\frac{y}{3} \le \abs{x} \le \frac{1}{2}} \abs{e(-(M+t) z) f_t(z) \prn{e\prn{-\tfrac{2t+1}{24} z}\frac{f_{t+1}(z)}{f_t(z)} - 1} dx}
\le 10 t y e^{-\frac{1}{70} \min(t, \frac{1}{y})} e^{2\pi M y} f_t(iy).
\qe

First, suppose $t y \le \frac{1}{2\pi} \log(4t)$. By Proposition \ref{prop:minor-arc-bound}, under the assumptions above on $t$ and $y$, we have
\aleq \label{eqn:diff-minor-arc-bound}
\int_{\frac{y}{3} \le \abs{x} \le \frac{1}{2}} \abs{f_t(z)} dx
&\le y e^{-\frac{1}{70} \min(t, \frac{1}{y})} f_t(iy)
\quad \text{and} \\
\int_{\frac{y}{3} \le \abs{x} \le \frac{1}{2}} \abs{f_{t+1}(z)} dx
&\le y e^{-\frac{1}{70} \min(t, \frac{1}{y})} f_{t+1}(iy).
\alqe

By Lemma \ref{lma:func_0-bound},
\eq
e^{2 \pi \frac{2t+1}{24} y} \frac{f_{t+1}(iy)}{f_t(iy)}
\le e^{7.5 t y e^{-2\pi t y}}
< 1.2,
\qe
which means,
\aleq
&\int_{\frac{y}{3} \le \abs{x} \le \frac{1}{2}} \abs{e(-(M+t) z) f_t(z) \prn{e\prn{-\tfrac{2t+1}{24} z}\frac{f_{t+1}(z)}{f_t(z)} - 1} dx} \\
&\quad \le 2.2 y e^{-\frac{1}{70} \min(t, \frac{1}{y})} e^{2 \pi (M+t) y} f_t(iy).
\alqe
Since $e^{2\pi t y} \le 4 t$, this proves (\ref{eqn:diff-away-from-0}) in this case.

Now, suppose $t y \ge \frac{1}{2\pi} \log(4t)$. By Lemma \ref{lma:func_0-bound}, we find
\eq
\abs{e\prn{-\tfrac{2t+1}{24} z}\frac{f_{t+1}(z)}{f_t(z)} - 1}
\le 7.5 \abs{t z} e^{-2\pi t y} e^{7.5 \abs{t z} e^{-2\pi t y}}
\le 10 t e^{-2\pi t y},
\qe
In this case, Proposition \ref{prop:minor-arc-bound} implies
\eq
\int_{\frac{y}{3} \le \abs{x} \le \frac{1}{2}} \abs{e(-(M+t) z) f_t(z) \prn{e\prn{-\tfrac{2t+1}{24} z} \frac{f_{t+1}(z)}{f_t(z)} - 1} dx}
\le 10 t y e^{-\frac{1}{70} \min(t, \frac{1}{y})} e^{2\pi M y} f_t(iy),
\qe
as in (\ref{eqn:diff-away-from-0}).

Our next goal is to replace $e\prn{-\tfrac{2t+1}{24} z}\frac{f_{t+1}(z)}{f_t(z)} - 1$ in the integrand with $-(2\pi i t z + 1) e(t z)$ in the range $\abs{x} < \frac{y}{3}$. To that end, let $\Upsilon(z) = \frac{f_{t+1}(z)}{f_t(z)} e\prn{-\tfrac{2t+1}{24} z}$, and apply Lemma \ref{lma:func_0-bound} to find
\aleq
& {\abs{\Upsilon(z) - 1 + e(tz) (2\pi i t z + 1)}} \\
&\le {\abs{\log\Upsilon(z) + e(tz) (2 \pi i t z + 1)}} + {\abs{\Upsilon(z) - 1 - \log\Upsilon(z)}} \\
&\le (40 \abs{z} + 22 e^{-2\pi t y}) \abs{t z} e^{-2\pi t y} + \frac{1}{2} \abs{7.5 t z e^{-2\pi t y}}^2 e^{7.5 \abs{tz} e^{-2\pi t y}} \\
&\le \prn{45 y + 85 t y e^{-2\pi t y}} t y e^{-2\pi t y}.
\alqe

Therefore, if we let
\eq
\a = \frac{\func_2(iy) - \func_2(ity)}{y} \in \prn{\frac{\min(t, \frac{1}{y})}{26}, \frac{\min(t, \frac{1}{y})}{12}}
\qe
as in \S\ref{ssec:asymptotic-proof}, then applying (\ref{eqn:f_t-abs-integral}), we have
\aleq \label{eqn:diff-simpler-integrand}
&\int_{-y/3}^{y/3} e(-(M+t) z) f_t(z) \prn{e\prn{-\tfrac{2t+1}{24} z} \frac{f_{t+1}(z)}{f_t(z)} - 1} dx \\
&= \int_{-y/3}^{y/3} e(-(M+t) z) f_t(z) e(tz) (-2\pi i t z - 1) dx \\
&\qquad + \err \sqrt{\frac{3}{2\a}} y e^{2\pi (M+t) y} f_t(iy) \prn{45 y + 85 ty e^{-2\pi t y}} t y e^{-2\pi t y} \\
&= \int_{-y/3}^{y/3} e(-M z) f_t(z) (-2\pi i t z - 1) dx \\
&\qquad + \err\prn{56 y + 105 t y e^{-2\pi t y}} \frac{t y^2}{\sqrt{\a}} e^{2\pi M y} f_t(iy).
\alqe

By Lemma \ref{lma:Gaussian-integral-with-error} and (\ref{eqn:reduction-to-gaussian}), we have
\aleq \label{eqn:diff-x-integral}
\abs{\int_{-y/3}^{y/3} e(-M z) f_t(z) (-2\pi i t x) dx}
&= 2\pi t y^2 e^{2\pi M y} f_t(iy) \abs{\int_{-1/3}^{1/3} x e^{-\pi \a x^2 \prn{1 + 2 i \varepsilon x + \err 3 x^2}} dx} \\
&\le 2 \pi t y^2 e^{2\pi M y} f_t(iy) \prn{2 \a^{-\frac{3}{2}}} \\
&\le 327 \frac{t y^3}{\sqrt{\a}} e^{2\pi M y} f_t(iy).
\alqe

What remains is
\aleq \label{eqn:diff-y-integral}
&\int_{-y/3}^{y/3} e(-M z) f_t(z) (2\pi t y - 1) dx \\
&= (2\pi t y - 1) \prn{\int_{-1/2}^{1/2} e(-M z) f_t(z)
+ \err 3 y e^{-\frac{1}{140}\frac{1}{y}} e^{2\pi M y} f_t(iy)} \\
&= (2\pi t y - 1) \prn{c_t(N)
+ \err y e^{-\frac{1}{70} \min(t, \frac{1}{y})} e^{2\pi M y} f_t(iy)}
\alqe

Altogether, combining (\ref{eqn:diff-away-from-0}, \ref{eqn:diff-simpler-integrand}, \ref{eqn:diff-x-integral}, \ref{eqn:diff-y-integral}), we have
\aleq \label{eqn:combined-diff-with-f_t}
c_{t+1}(N+t) - c_t(N+t)
&= (2\pi t y - 1) c_t(N) + \err t y e^{2\pi M y} f_t(iy) \\
&\qquad \cdot \prn{
11 e^{-\frac{1}{70} \min(t, \frac{1}{y})}
+ 383 \frac{y^2}{\sqrt{\a}}
+ 105 \frac{t y^2}{\sqrt{\a}} e^{-2\pi t y}
}.
\alqe

Finally, note that by Theorem \ref{thm:t-core-asymptotic-explicit},
\eq
c_t(N)
= \frac{y}{\sqrt{\a}} e^{2\pi M y} f_t(iy) \prn{1 + \err \frac{3.51}{\a}}
> 0.9 \frac{y}{\sqrt{\a}} e^{2\pi M y} f_t(iy),
\qe
so we may write (\ref{eqn:combined-diff-with-f_t}) as
\aleq
&c_{t+1}(N+t) - c_t(N+t) \\
&\quad = c_t(N) \prn{
2\pi t y - 1 + \err \frac{t y}{0.9} \prn{
    11 \frac{\sqrt{\a}}{y} e^{-\frac{1}{70} \min(t, \frac{1}{y})}
    + 383 y
    + 105 t y e^{-2\pi t y}
    }
} \\
&\quad = c_t(N) \prn{
2 \pi t y - 1 + \err t y \prn{705 y + 120 t y e^{-2\pi t y}}},
\alqe
as claimed.
\end{proof}

\subsection{Proof of Theorem \ref{thm:t-core-asymptotic-small-t}} \label{ssec:small-t-proof}

As in \S\ref{ssec:asymptotic-proof}, we may assume $N$ and $t$ are sufficiently large. Let $\hat{\e} = \min(\e, 3)$ and choose $y = \frac{t-1}{4\pi M}$, which we may assume is sufficiently small. The hypotheses in Theorem \ref{thm:t-core-asymptotic-small-t} imply
\eq \label{eqn:small-t-hypothesis}
t y < \frac{2\pi}{(2+\hat{\e}) \log t}.
\qe
We may repeat the proof of Proposition \ref{prop:minor-arc-bound} to find
\eq
\int_{y \sqrt{1+\hat{\e}} \le \abs{x} \le \frac{1}{2}} \abs{\frac{f_t(z)}{f_t(iy)}} dx
\ll \frac{y}{\sqrt{t-1}} e^{-\frac{1}{6} t}
\qe

The hypothesis (\ref{eqn:small-t-hypothesis}) also implies that if $\abs{x} \le y \sqrt{1+\hat{\e}}$, then
\eq
\Im \frac{-1}{tz} \ge \frac{1}{(2+\hat{\e}) t y} > \frac{\log t}{2\pi},
\qe
which means
\aleq
f_t(z)
&= t^{-\frac{t}{2}} (-iz)^{\frac{1-t}{2}} \frac{\eta\prn{\frac{-1}{tz}}^t}{\eta\prn{\frac{-1}{z}}} \\
&= t^{-\frac{t}{2}} (-iz)^{\frac{1-t}{2}} e^{\err 2.02 t e^{-2\pi \Im \frac{-1}{tz}}} \\
&= t^{-\frac{t}{2}} (-iz)^{\frac{1-t}{2}} \prn{1 + O\prn{t e^{-2\pi \Im \frac{-1}{tz}}}},
\alqe
and in particular $f_t(iy) \le 2 t^{-\frac{t}{2}} y^{\frac{1-t}{2}}$.

Our goal is to evaluate
\aleq \label{eqn:small-t-integrals}
c_t(N)
&= t^{-\frac{t}{2}} \int_{-y\sqrt{1+\hat{\e}}}^{y\sqrt{1+\hat{\e}}} e(-M z) (-iz)^{\frac{1-t}{2}} dx \\
&\quad + O\prn{
    t^{\frac{2-t}{2}} \int_{-y\sqrt{1+\hat{\e}}}^{y\sqrt{1+\hat{\e}}} \abs{e(-M z) (-iz)^{\frac{1-t}{2}} e^{-2\pi \Im \frac{-1}{tz}}} dx
    } \\
&\quad + O\prn{
    e^{2\pi M y} t^{-\frac{t}{2}} y^{\frac{1-t}{2}} \frac{y}{\sqrt{t-1}} e^{-\frac{1}{6} t}
    }.
\alqe

Let $\mathcal{L}$ be the contour that traverses in straight lines from $-y\sqrt{1+\hat{\e}} - i\infty$ to $-y\sqrt{1+\hat{\e}} + iy$ to $y\sqrt{1+\hat{\e}} + iy$ to $y\sqrt{1+\hat{\e}} - i \infty$. Hankel's formula states
\eq
\int_{\mathcal{L}} e(-M z) (-i z)^{\frac{1-t}{2}} d z = \frac{(2\pi)^{\frac{t-1}{2}}}{\G\prn{\frac{t-1}{2}}} M^{\frac{t-3}{2}},
\qe
which means that the first line in the right-hand side of (\ref{eqn:small-t-integrals}) is
\eq \label{eqn:small-t-integrals-bound-1}
\frac{(2\pi)^{\frac{t-1}{2}}}{t^{\frac{t}{2}} \G\prn{\frac{t-1}{2}}} M^{\frac{t-3}{2}} + O\prn{
    e^{2\pi M y} t^{-\frac{t}{2}} y^{\frac{1-t}{2}} \frac{(1+\hat{\e})^{\frac{1-t}{4}}}{M}
    }.
\qe
By widening the domain of integration from $\abs{x} \le y \sqrt{1+\hat{\e}}$ to $\abs{x} \le 2 y$ and making the change of variables $u = \frac{x}{y}$, we find that the second line in the right-hand side of (\ref{eqn:small-t-integrals}) is
\eq \label{eqn:small-t-integrals-bound-2}
\ll e^{2\pi M y} t^{\frac{2-t}{2}} y^{\frac{3-t}{2}} \int_{-2}^{2} (1+u^2)^{\frac{1-t}{4}} e^{-\frac{2\pi}{ty} \frac{1}{1+u^2}} du.
\qe
We will now show
\eq \label{eqn:small-t-error-integral}
\int_{-2}^{2} (1+u^2)^{\frac{1-t}{4}} e^{-\frac{2\pi}{ty} \frac{1}{1+u^2}} du
\ll \max\prn{e^{-\frac{1}{5} t}, \frac{1}{\sqrt{t-1}} e^{-\frac{2\pi}{ty}}}.
\qe
Indeed, the function $v \mapsto v^{\frac{t-1}{4}} e^{-\frac{2\pi}{ty} v}$ is maximized at $v = \frac{t(t-1)y}{8\pi}$. Therefore, if $\frac{t(t-1)y}{8\pi} \le 1$, then the integrand in (\ref{eqn:small-t-error-integral}) is bounded by $e^{-\frac{1}{4} t}$, while if $\frac{t(t-1)y}{8\pi} \ge 1$, then the integrand in (\ref{eqn:small-t-error-integral}) is bounded by $e^{-\frac{2\pi}{ty}}$ and decreasing in $u$. This is sufficient to prove (\ref{eqn:small-t-error-integral}) when $\frac{t(t-1)y}{8\pi} \le \frac{5}{4}$. When $\frac{t(t-1)y}{8\pi} \ge \frac{5}{4}$, note that $\frac{1}{1+u^2} \ge 1 - u^2$ and $1 + u^2 \ge e^{\frac{5}{6} u^2}$ if $\abs{u} \le \frac{1}{2}$. Therefore,
\aleq
\int_{-2}^{2} (1+u^2)^{\frac{1-t}{4}} e^{-\frac{2\pi}{ty} \frac{1}{1+u^2}} du
&\ll \int_{-1/2}^{1/2} e^{-\frac{5}{24} (t-1) u^2} e^{-\frac{2\pi}{ty} (1 - u^2)} du \\
&\le e^{-\frac{2\pi}{ty}} \int_{-\infty}^\infty e^{-\frac{1}{120} (t-1) u^2} du \\
&\ll \frac{1}{\sqrt{t-1}} e^{-\frac{2\pi}{ty}},
\alqe
proving (\ref{eqn:small-t-error-integral}).

Applying (\ref{eqn:small-t-integrals-bound-1}, \ref{eqn:small-t-integrals-bound-2}, \ref{eqn:small-t-error-integral}) to (\ref{eqn:small-t-integrals}), we find
\aleq
c_t(N)
&= \frac{(2\pi)^{\frac{t-1}{2}}}{t^{\frac{t}{2}} \G\prn{\frac{t-1}{2}}} M^{\frac{t-3}{2}} \\
&\quad + O\prn{e^{2\pi M y} t^{-\frac{t}{2}} y^{\frac{1-t}{2}} \prn{
\frac{(1+\hat{\e})^{\frac{1-t}{4}}}{M}
+ t y \max\prn{e^{-\frac{1}{5} t}, \frac{1}{\sqrt{t-1}} e^{-\frac{2\pi}{ty}}}
+ \frac{y}{\sqrt{t-1}} e^{-\frac{1}{6} t}
}}.
\alqe

Recalling the choice $y = \frac{t-1}{4\pi M}$ (which minimizes $e^{2\pi M y} y^{\frac{1-t}{2}}$) and applying Stirling's formula, we find
\eq
e^{2\pi M y} t^{-\frac{t}{2}} y^{\frac{1-t}{2}}
= t^{-\frac{t}{2}} \prn{\frac{4e\pi M}{t-1}}^{\frac{t-1}{2}}
\ll \frac{\sqrt{t-1}}{y} \frac{(2\pi)^{\frac{t-1}{2}}}{t^{\frac{t}{2}} \G\prn{\frac{t-1}{2}}} M^{\frac{t-3}{2}},
\qe
and hence
\eq
c_t(N)
= \frac{(2\pi)^{\frac{t-1}{2}}}{t^{\frac{t}{2}} \G\prn{\frac{t-1}{2}}} M^{\frac{t-3}{2}}
\prn{1 + O\prn{
    t^{-\frac{1}{2}} (1+\e)^{-\frac{1}{4} t}
    + e^{-\frac{1}{6} t}
    + t e^{-\frac{2\pi}{ty}}
    }
},
\qe
as claimed.

\begin{RMK} \label{rmk:t-core-asymptotic-small-t-explicit}
In \S\ref{sec:stanton}, we will require an explicit version of Theorem \ref{thm:t-core-asymptotic-small-t} to prove Stanton's conjecture for $8 \le t \le 1000$. By keeping careful track of constants and repeating the proof above with $\e = 3$, we can show that if $t \ge 8$, $M \ge 100\,000$, and $\frac{t(t-1)}{4 \pi M} < \frac{1}{5} \min\prn{\frac{2}{\sqrt{3}}, \frac{2\pi}{\log t}}$, then
\eq
c_t(N)
= \frac{(2\pi)^{\frac{t-1}{2}}}{t^{\frac{t}{2}} \G\prn{\frac{t-1}{2}}} M^{\frac{t-3}{2}} \prn{1 + \err \prn{2.5 e^{-\frac{1}{8} t} + 20 t^{-4}}}.
\qe
\end{RMK}

\subsection{Proof of Theorem \ref{thm:t-core-asymptotic-big-t}} \label{ssec:big-t-proof}

We may assume $N \ge 3t$, since otherwise $c_t(N) = p(N) - t p(N-t) + \frac{1}{2} (t^2-3t) p(N-2t)$.

Let
\eq
y = \frac{1}{\frac{3}{\pi} + \sqrt{24(N-2t) - 1 + \frac{9}{\pi^2}}},
\qe
so that our assumption $t > (1+\e) \frac{\sqrt{6}}{2\pi} \sqrt{N} \log N$ implies $\abs{t e(tz)} = t e^{-2\pi t y} \ll_\e N^{-\frac{\e}{2}} \log N$.

Consider the integrals
\aleq
p(N)
&= \int_{-1/2}^{1/2} e\prn{-\prn{N-\tfrac{1}{24}}z} \frac{1}{\eta(z)} dx
\quad \text{and} \\
c_t(N)
&= \int_{-1/2}^{1/2} e\prn{-\prn{N-\tfrac{1}{24}}z} \frac{1}{\eta(z)} e\prn{-\tfrac{t^2}{24} z}\eta(tz)^t dx.
\alqe
Note that
\eq
\log \prn{e\prn{-\tfrac{t^2}{24} z}\eta(tz)^t}
= -\frac{2\pi i}{z} \prn{\func_0(tz) + \frac{(tz)^2}{24}}
= -t e(tz) + O_\e\prn{t e(2tz)},
\qe
so
\aleq
c_t(N)
&= \int_{-1/2}^{1/2} e\prn{-\prn{N-\tfrac{1}{24}}z} \frac{1}{\eta(z)} e^{-t e(tz) + O_\e\prn{t e(2tz)}} dx \\
&= \int_{-1/2}^{1/2} e\prn{-\prn{N-\tfrac{1}{24}}z} \frac{1}{\eta(z)} \prn{1 - t e(tz) + O_\e\prn{t^2 e(2tz)}} dx \\
&= p(N) - t p(N-t) + O_\e\prn{t^2 \int_{-1/2}^{1/2} \abs{e\prn{-\prn{N-2t-\tfrac{1}{24}}z} \frac{1}{\eta(z)}} dx} \\
&= p(N) - t p(N-t) + O_\e\prn{t^2 p(N-2t)},
\alqe
with this last equality implied by e.g. (\ref{eqn:f_t-abs-integral}) applied to $c_N(N-2t) = p(N-2t)$.

\section{Proof of Stanton's conjecture}
\label{sec:stanton}
In this section, we will prove Stanton's conjecture \ref{conj:stanton}.

Recall from the introduction that we may assume $4 \le t < N-1$. To begin, we check using SageMath \cite{sagemath} that Stanton's conjecture holds for all $t < N - 1$ and $N < 100\,000$, and also that it holds for all $t < 1000$ and $N < 500\,000$.

Kim-Rouse \cite{kim-rouse} verified the conjecture for $4 \le t < 8$. In \S\ref{ssec:stanton-small-t}, we will use Theorem \ref{thm:t-core-asymptotic-small-t} (or rather, the explicit version mentioned in Remark \ref{rmk:t-core-asymptotic-small-t-explicit}) to verify the conjecture for $8 \le t < 1000$ and $N \ge 500\,000$.

If $t \ge 1000$ and $N \ge 100\,000$, then as in \S\ref{ssec:y-def}, choose $y$ in the range
\eq
\frac{t-1}{4\pi M}
< y <
\frac{1}{\frac{3}{\pi} + \sqrt{24N-1 + \frac{9}{\pi^2}}} < \frac{1}{1000}
\qe
such that
\eq
\frac{\func_1(ity) - \func_1(iy)}{y^2} = M.
\qe
We will complete the proof of Stanton's conjecture by treating the case $t y \le \frac{1}{2}$ in \S\ref{ssec:stanton-medium-t} using Theorem \ref{thm:t-core-asymptotic-explicit}, and the case $t y \ge \frac{1}{2}$ in \S\ref{ssec:stanton-big-t} using Theorem \ref{thm:t-core-diff-asymptotic-explicit}.

\begin{RMK}
For any $\e > 0$ and $t$ and $N$ sufficiently large depending on $\e$, we could prove Stanton's conjecture using Theorem \ref{thm:t-core-asymptotic-small-t} if $t y < (1 - \e) \frac{\pi}{\log t}$, using Theorem \ref{thm:t-core-asymptotic-explicit} if $t y < (1 - \e) \frac{\log t}{2\pi}$, using Theorem \ref{thm:t-core-diff-asymptotic-explicit} if $t y \ge \frac{1}{2}$, and using Theorem \ref{thm:t-core-asymptotic-big-t} if $t y > (1 + \e) \frac{\log t}{\pi}$, giving us some leeway in deciding which theorem to use for which range. Our decision was based on an attempt to minimize the number of cases that must be checked by hand.
\end{RMK}

\subsection{Stanton's conjecture for \texorpdfstring{$8 \le t < 1000$}{8 <= t < 1000} and \texorpdfstring{$N \ge 500\,000$}{N >= 500000}} \label{ssec:stanton-small-t}

With these assumptions on $t$ and $N$, we have $\frac{t(t-1)}{4 \pi M} < \frac{1}{5} \min\prn{\frac{2}{\sqrt{3}}, \frac{2\pi}{\log t}}$, so by Remark \ref{rmk:t-core-asymptotic-small-t-explicit},
\eq
c_t(N)
= \frac{(2\pi)^{\frac{t-1}{2}}}{t^{\frac{t}{2}} \G\prn{\frac{t-1}{2}}} M^{\frac{t-3}{2}} \prn{1 + \err E(t)}
\qe
where $E(t) = 2.5 e^{-\frac{1}{8} t} + 20 t^{-4}$.

Applying the inequalities
\eq
\frac{t^{\frac{t}{2}}}{(t+1)^{\frac{t+1}{2}}} > \frac{1}{\sqrt{e(t+1)}}
\quad \text{and} \quad
\frac{\G\prn{\frac{t-1}{2}}}{\G\prn{\frac{t}{2}}} > \sqrt{\frac{2}{t}}
\qe
yields
\eq
\frac{c_{t+1}(N)}{c_t(N)}
> \sqrt{\frac{4\pi M}{e t (t+1)}} \frac{1 - E(t+1)}{1 + E(t)}.
\qe
For all $8 \le t < 1000$, we have
\eq
\sqrt{\frac{4\pi \cdot 500\,000}{e t (t+1)}} \frac{1 - E(t+1)}{1 + E(t)}
> 1.5,
\qe
which proves Stanton's conjecture in this range.

\subsection{Stanton's conjecture for \texorpdfstring{$t \ge 1000$}{t >= 1000}, \texorpdfstring{$N \ge 100\,000$}{N >= 100000}, and \texorpdfstring{$t y \le \frac{1}{2}$}{ty <= 1/2}} \label{ssec:stanton-medium-t}

Under these assumptions, we may use Theorem \ref{thm:t-core-asymptotic-explicit}, which tells us
\eq
c_t(N) = \frac{y^{\frac{3}{2}} e^{2 \pi M y} f_t(iy)}{\sqrt{\func_2(iy) - \func_2(ity)}}
\prn{
1 + \err \frac{3.5 y}{\func_2(iy) - \func_2(ity)}
}.
\qe
This same value of $y$ also satisfies the hypotheses in Theorem \ref{thm:t-core-asymptotic-explicit} for $c_{t+1}(N)$. Indeed, by the mean value theorem,
\aleq
\abs{y \prn{\frac{\func_1(i(t+1)y) - \func_1(iy)}{y^2} - \prn{N + \frac{(t+1)^2-1}{24}}}}
&= \abs{\frac{\func_1(i(t+1) y) - \func_1(ity)}{y} - \frac{2t+1}{24} y} \\
&= \abs{i \func_1'(isy) - \frac{sy}{12}}
\alqe
for some $s \in (t, t+1)$, and by (\ref{eqn:func-prime-fourier}),
\aleq 
\abs{i \func_1'(isy) - \frac{sy}{12}}
= \abs{\frac{1}{4\pi} - \frac{sy}{12} + \sum_{n=1}^\infty \prn{\frac{2\pi n}{sy}}^2 \frac{\s(n)}{2\pi n} e^{-2\pi n/(sy)}}
< \frac{1}{4\pi}.
\alqe
Therefore, Theorem \ref{thm:t-core-asymptotic-explicit} also tells us
\eq
c_{t+1}(N) = \frac{y^{\frac{3}{2}} e^{2 \pi \prn{M + \frac{2t+1}{24}} y} f_{t+1}(iy)}{\sqrt{\func_2(iy) - \func_2(i(t+1)y)}}
\prn{
1 + \err \frac{3.5 y}{\func_2(iy) - \func_2(i(t+1)y)}
}.
\qe

By part \ref{item:func_2-prime-increasing} of Lemma \ref{lma:func_2-bound}, we have
\eq
\frac{3.5 y}{\func_2(iy) - \func_2(ity)}
\le \frac{3.5}{-i (t-1) \func_2'(\frac{i}{2})}
\le \frac{45}{t}
\qe
and
\eq
\frac{3.5 y}{\func_2(iy) - \func_2(i(t+1)y)}
\le \frac{3.5}{-i t \func_2'(\frac{501i}{1000})}
\le \frac{45}{t},
\qe

so
\eq \label{eqn:stanton-ratio-1}
\frac{c_{t+1}(N)}{c_t(N)}
= e^{2 \pi y \frac{2t+1}{24}} \frac{f_{t+1}(iy)}{f_t(iy)} \sqrt{\frac{\func_2(iy) - \func_2(ity)}{\func_2(iy) - \func_2(i(t+1)y)}} \frac{1 + \err \frac{45}{t}}{1 + \err \frac{45}{t}}.
\qe

As in (\ref{eqn:func_0-mvt}), there is some $s \in (t, t+1)$ such that
\aleq \label{eqn:stanton-ratio-2}
\log\prn{e^{2 \pi \prn{\frac{2t+1}{24}} y} \frac{f_{t+1}(iy)}{f_t(iy)}}
&= -2\pi i \prn{\func_0'(isy) + \frac{isy}{12}} \\
&= -\frac{1}{2} \log(e s y) + \frac{\pi}{6} s y - \sum_{n=1}^\infty \prn{\frac{2\pi n}{s y} + 1} \frac{\s(n)}{n} e^{-2\pi n /(sy)} \\
&> 0.1.
\alqe
We also have
\eq \label{eqn:stanton-ratio-3}
\frac{\func_2(iy) - \func_2(ity)}{\func_2(iy) - \func_2(i(t+1)y)}
= 1 - \frac{\func_2(ity) - \func_2(i(t+1)y)}{\func_2(iy) - \func_2(i(t+1)y)}
= 1 + \err \frac{1.1}{t},
\qe

so altogether (\ref{eqn:stanton-ratio-1}, \ref{eqn:stanton-ratio-2}, \ref{eqn:stanton-ratio-3}) imply
\eq
\frac{c_{t+1}(N)}{c_t(N)}
> e^{0.1} \sqrt{1 - \frac{1.1}{1000}} \frac{1 - \frac{45}{1000}}{1 + \frac{45}{1000}}
> 1.
\qe

\subsection{Stanton's conjecture for \texorpdfstring{$t \ge 1000$}{t >= 1000}, \texorpdfstring{$N \ge 100\,000$}{N >= 100000}, and \texorpdfstring{$t y \ge \frac{1}{2}$}{t y >= 1/2}} \label{ssec:stanton-big-t}
With these assumptions, let us choose $y'$ such that
\eq
\frac{\func_1(ity')-\func_1(iy')}{y'^2} = N - t + \frac{t^2-1}{24}.
\qe
Since $\frac{\func_1(ity)-\func_1(iy)}{y^2}$ is monotonically decreasing in $y$ by (\ref{eqn:y-unique}), we have $y' > y \ge \frac{1}{2t}$.

If $t > \frac{N}{2}$ (so $c_t(N) = p(N) - t p(N-t)$), then Stanton's conjecture was proved by Craven \cite{craven}. Otherwise, $N - t \ge \frac{N}{2} \ge 50\,000$, so
\eq
\frac{1}{y'} > \frac{3}{\pi} + \sqrt{24\cdot50\,000 - 1 + \frac{9}{\pi^2}} > 1000,
\qe
and Theorem \ref{thm:t-core-diff-asymptotic-explicit} applies, which means
\aleq
c_{t+1}(N) - c_t(N)
&= c_t(N-t) \prn{2\pi t y' - 1 + \err t y' \prn{705 y' + 120 t y' e^{-2\pi t y'}}}, \\
&\ge c_t(N-t) \prn{2\pi t y' - 1 - t y' \prn{\frac{705}{1000} + 60 e^{-\pi}}} \\
&> 0.
\alqe

\pagebreak
\printbibliography


\end{document}